\documentclass[a4]{article}
\usepackage[english]{babel}
\usepackage{amsmath,amssymb,amsfonts,latexsym,wasysym}
\usepackage{stmaryrd}
\usepackage{array}
\usepackage{ulem}
\usepackage{amsthm}
\usepackage{fancybox,pifont,calc}
\usepackage[latin1]{inputenc}
\usepackage[T1]{fontenc}
\usepackage{graphicx}

\newcommand{\oo}{\infty}

\newcommand{\Z}{\mathbb{Z}}

\newcommand{\N}{\mathbb{N}}

\newcommand{\beq}{\begin{equation}}
\newcommand{\eeq}{\end{equation}}
\newcommand{\beqn}{\begin{eqnarray}}
\newcommand{\beqnn}{\begin{eqnarray*}}
\newcommand{\eeqn}{\end{eqnarray}}
\newcommand{\eeqnn}{\end{eqnarray*}}

\newtheorem{thm}{Theorem}[section]

\newtheorem{lem}[thm]{Lemma}

\begin{document}

\title{ A new enlightenment about gaps between primes} 
\date{June 16, 2015}
\author{Belhaouari, Samir B.  \\
	Department of Computer science, INNOPOLIS University, Kazan, Russia\\
	Department of Mathematics, University of Sharjah, UAE \\
email: samir.brahim@gmail.com\\}
	
\maketitle

\begin{abstract}
The idea of generating prime numbers through sequence of sets of co-primes was the starting point of this paper that ends up by proving two conjectures, the existence of infinitely many twin primes and the Goldbach conjecture.
The main idea of our approach is summarized on the creation and on the analyzing sequence of sets of distinct co-primes with the first $n$ primes, $\left\{ p_i :\, i\leq n \right\}$, and the important properties of the modulus linear combination of the co-prime sets, $H=\left(1,p_{n+1},..., \Pi_{i=1}^n p_i-1\right) $, that gives sets of even numbers $\{0,2,4,..., \Pi_{i=1}^n p_i -2 \}$.
Furthermore, by generalizing our approach, the Polignac conjecture "the existence of infinitely many cousin primes, $p_{n+1}-p_{n}=4$," and the statement that "every even integer can be expressed as a difference of two primes," are derived as well.
\end{abstract}



\tableofcontents

\section{Introduction}
    The proof of the existence of infinitely many twin primes has been one of the most difficult and the oldest question in number theory for several years.
The twin prime conjecture states: There are infinitely many primes $p$ such that $p + 2$ is also prime. Some attribute the conjecture to the Greek mathematician Euclid of Alexandria, who gave the oldest known proof that there exist an infinite number of primes.\\
In 1849 de Polignac made the more general conjecture that for every integer $k$, there are infinitely many prime pairs $(p_1, p_2) $ such that $p_1 - p_2 = 2k$. The twin prime conjecture is for $k = 1$.

Since every three consecutive numbers, $\left(p, p+1, p+2 \right)$, are not all co-primes with the number 3, then 2 and 3 are factors of $p+1$. Subsequently, all twin primes are in the form of $6n\pm1$ except for the first pair $(3,5)$ and also all cousin primes are in the form of $6n+1$ and $6n+5$ consecutively. On December 2011, in \cite{paper3}, the largest known twin prime pair, $3756801695685*2^{666669} \pm 1$, were found with 200700 digits.

Similar to the distribution of prime numbers among positive integers, the number of twin primes less that $N$ naturally declines as $N$ increases. This brings up the question of how many pairs of twin primes are there?

Many research works have been devoted to provide a proof to the twin prime conjecture but it is still one of the open problems in mathematics and number theory.
Some important result of twin primes and the gap between primes have been obtained. One of these results is Brun's constant (B). In 1919 Brun, \cite{Brun}, proved that the sum of the reciprocals of the twin primes converge to a definite number called Brun's constant and has the value of $B = 1.90218 \pm 2 \cdot 10^{-5}$.
The recently breakthrough, in 2013, of Yitang Zhang \cite{paper4}, where he proved that there are infinitely many pairs of primes with prime gap less than 70 million
$$\liminf_{n\rightarrow \oo} \left(p_{n+1}-p_n \right) <7.10^{7}.$$

This proof shows that the gap between primes does not grow infinitely, and it is bounded by a definite number.
This result was significantly improved by several researchers.\
\newline

In May 2013, mathematicians had uncovered simple tweaks to Zhang's argument that brought the bound below 60 million. A May 30 blog post by Scott Morrison of the Australian National University in Canberra ignited a firestorm of activity, as mathematicians vied to improve on this number, setting one record after another. By June 4, Terence Tao of the University of California, Los Angeles, a winner of the Fields Medal, mathematics' highest honor, had created a Polymath project, an open, online collaboration to improve the bound that attracted dozens of participants.
For weeks, the project moved forward at a breathless pace. At times, the bound was going down every thirty minutes, Tao recalled. By July 27, the team had succeeded in reducing the proven bound on prime gaps from 70 million to 4,680.

On November 19 by James Maynard \cite{Maynard}, a postdoctoral researcher working on his own at the University of Montreal, has upped the ante. Just months after Zhang announced his result, Maynard has presented an independent proof that pushes the gap down to 600. A new Polymath project is in the planning stages, to try to combine the collaboration's techniques with Maynard's approach to push this bound even lower, for more details please visit the website WIRED, Nov. 2013.

Strong connection between the twin prime and Goldbach's Conjectures was found along of our approach to prove the twin prime's conjecture. The Goldbach conjecture, dating from Goldbach's correspondence with Euler in 1742, that states: Every even integer greater than 2 is the sum of two prime numbers (not
necessarily distinct). It has been verified that the even integers through $10^8$ have the stated property, but Goldbach's conjecture remains unproved.

Several works on background information about Goldbach's conjecture have been conducted previously. Namely Dickson \cite{Dickson} lists results obtained through the first decade of this century. Wang, \cite{Wang}, is a collection of this century's most important papers on the Goldbach conjecture through the early 1980's. Shanks, \cite{Daniel}, is a marvelous book on conjectures in number theory, Goldbach's among them, and on mathematical conjectures in general. Ribenboim, \cite{Paulo}, is another treasure chest in number theory, with a section on Goldbach (pp. 229-235).\
\newline
\\
\section{Abbreviation and Notation}
Before starting the proof, some important definitions are needed to be introduced as it is indicated below 

\begin{description}
\item[$\centerdot$]  Let $p_2$,..., $p_n$ distinct ordered prime numbers, and for only notation reason we put $p_1=1$.
  \item[$\centerdot$] For any two sets $U$ and $V$ of integer numbers, their sum and their difference is  defined as follows
                 $$ U\pm V= \{x\pm y \;\; : x\in U, y\in V\}.$$
  \item[$\centerdot$] For any set $U$ is multiplied by a value $c$ is defined as follows $$c\,U=\{c.x\; :\; x\in U\}.$$
  \item[$\centerdot$] The difference between all elements of the set $U$, is noted by $$\Delta U=\{x-y\; : (x,y)\in U^2\}.$$
   \item[$\centerdot$] The set $U/V$ is set of elements in $U$ but not in $V$.
  \item[$\centerdot$] The arithmetic sequence, with reason $r$, will be noted by $[r_0:r:r_f]$, where $r_0$ and $r_f$ are the initial and the final values of the sequence respectively, i.e.,
      $$[r_0:r:r_f]=\{r_0, r_0+r, r_0+2r,...,r_f\}.$$
  \item[$\centerdot$] If the the reason of the arithmetic sequence is equal to one then it will be noted simply by $[r_0:r_f]$, i.e.,
                             $$[r_0:r_f]=\{r_0, r_0+1, r_0+2,...,r_f\}.$$
  \item[$\centerdot$] The function $\lfloor x \rfloor$ gives the largest integer less than or equal to x.
  \item[$\centerdot$] The modulo operation is computing the remainder after division, $m$ Modulo $n$ will be noted  as $m_{mod(n)}$, i.e.,
                $$m_{mod(n)}=m-n \left\lfloor\frac{m}{n}\right\rfloor.$$
  \item[$\centerdot$] The set $V_{x}$ contains all co-prime numbers with all devisors of the integer $x$ and are less than $x$.
  \item[$\centerdot$] The set $V_{\prod_{i=1}^{n} p_i}$ contains all co-prime numbers with primes $p_i$, $i=\{1,2,...,n\}$, and less than $\prod_{i=1}^{n} p_i$.
  \item[$\centerdot$] The set $w_k$ is the set of co-primes with  primes $p_i, \; i=\{n-k+2,,...,n+1\}$, and less then $\prod_{i=n-k+2}^{n+1} p_i$.
  \item[$\centerdot$] Let $\pi(x)$ be the prime-counting function that gives the number of primes less or equal to $x$.
  \item[$\centerdot$] Let $\pi(x,n)$ be the co-prime-counting function that gives the number of co-primes with $p_i, i>n$, and are less or equal to $x$.
	\item [$\centerdot$] Let $\Psi(x,y,Q)$ be the number of all positive integers less than $x$ which do not have a prime divisor which is less than $y$ and belongs to $Q$.
	
  \item[$\centerdot$] $U^k_{\prod_{i=1}^{n+1} p_i}=V_{\prod_{i=1}^{n+1} p_i} \bigcap  \left[0:  \prod_{i=1}^{n+1-k} p_i \right] $, unless it is defined differently.
    \item[$\centerdot$] $U^{k,j}_{\prod_{i=1}^{n+1} p_i}=V_{\prod_{i=1}^{n+1} p_i} \bigcap  \left[j \prod_{i=1}^{n+1-k} p_i: (j+1) \prod_{i=1}^{n+1-k} p_i \right] $.
  \item[$\centerdot$] The cardinality of a set $U$ is noted by $Card (U).$
  \item[$\centerdot$] If the asymptomatic ration of two functions verify
	 $$\lim_{n \rightarrow \infty} \left|\frac{f(n)}{g(n)}\right|=1,$$ then it will be noted by
  $$f\sim g.$$ 
	\item [$\centerdot$] If there exist positive $\delta$ and $M$ such that
	$$\left|f(x)\right| \leq M \left| g(x) \right| \; \text{ for } x> \delta,$$
	then we write
	$$f(x)=O(g(x)) \;\text{ as } x\rightarrow \infty.$$
	\item[$\centerdot$] The greatest common divisor pf two integers, n and m, will be noted by
	$$gcd(n,m).$$
\end{description}

\section{Sketch of the proof}
This section will summarized all the headlines of our approach, from the starting point, through important properties, till reaching the two conjectures. All headlines are briefed as follows
\begin{itemize}
  \item Any odd number $n$ can be written under the form of $n=2k+1$, where $k \in \N$.
	
  \item All co-primes numbers with 2 and 3 can be written under the form $6k+1$ or $6k+5$, so the two forms can be merged under one form $n=6k+V_{2*3}$, where $k \in \N$ and the set $V_{2*3}=\{1,5\}$.
	
  \item All co-primes number with 2, 3, and 5 can be written under the form of $n=30k+V_{2*3*5}$, where $k \in \N$ and $V_{2*3*5}=\{1,7,11,13,17,19,23,29\}$.
	
  \item Theorem \ref{theorem 1} generalizes this idea of finding the shape of all co-primes number with $p_i, i\leq n+1$, as
  $$n= k \left(\Pi_{i=1}^{n+1} p_i\right)+ V_{mod\left ( \prod_{i=1}^{n+1}P_{i} \right )},$$
  where the sets $V_.$ are generated by the following recurrence formula
  $$V_{mod\left ( \prod_{i=1}^{n+1}P_{i} \right )}=\left(P_{n+1}V_{mod\left ( \prod_{i=1}^{n}P_{i} \right )}+\left ( \prod_{i=1}^{n}P_{i} \right )\left [ 1:P_{n+1}-1 \right ]\right)_{mod\left ( \prod_{i=1}^{n+1}P_{i} \right )}.$$
  
	\item We have noticed and proved that for all $n\geq2$
  $$\left(V_{\prod_{i=1}^{n} p_i} - V_{\prod_{i=1}^{n} p_i}\right)_{mod\left(\prod_{i=1}^{n} p_i\right)}=\left[0:2:\prod_{i=1}^{n} p_i-2 \right].$$
  
	\item We have noticed also that $\forall\, n\in \N \;
    $ \begin{eqnarray*}
                        V_{\prod_{i=1}^{n} p_i} &=&\left\{1,\,p_{n+1},..., \left(\prod_{i=1}^{n} p_i\right) -p_{n+1}, \, \left(\prod_{i=1}^{n} p_i\right)-1 \right\} \\
                         \min V_{\prod_{i=1}^{n} p_i}/\{1\} &=& p_{n+1} \\
                          V_{\prod_{i=1}^{n} p_i} \bigcap (1: p_{n+1}) &=& \emptyset.
                       \end{eqnarray*}
											
  \item The set $ V_{\prod_{i=1}^{n} p_i} \cap (0: p_{n+1}^2)$ contains only primes numbers.
  
	\item The asymptotic density of co-primes are found to be as follows
  \begin{eqnarray*}
    \pi(x,n) &\sim& \ln(n) \frac{x}{\ln(x)}, \\
    \frac{\pi(x,n)}{x} &\sim&  \frac{\ln(n)}{\ln(x)} ,\\
    &\text{and}&  \\
    \frac{\text{card}(V_{\prod_{i=1}^{n} p_i})}{\prod_{i=1}^{n} p_i}  &=&  \frac{\Pi_{i=1}^n (p_i-1)}{\Pi_{i=1}^n p_i}\\
    &\sim & \exp\left(\int^{n} \frac{-dn}{n\ln(n)}\right)\\
    &\sim & \frac{C}{\ln(n)}\\
    &\sim & C\frac{\pi(n)}{n}.\\
  \end{eqnarray*}

  \item In order to investigate the difference between only primes number, the primes inside the $V_{\prod_{i=1}^{n} p_i}$ will be considered as it is indicated in Lemma \ref{LemmaP}, which it gives extremely important properties as follows:
  \begin{equation*}
\left(\Delta V_{\prod_{i=1}^{n+1} p_i}\right)_{mod(\prod_{i=1}^{n+1} p_i)}=\left(\Delta U^k_{\prod_{i=1}^{n+1} p_i}\right)_{mod(\prod_{i=1}^{n+1-k} p_i)}+\sum_{j=0}^{k-1} \left(\prod_{i=1}^{n-j} p_i \right)[0:p_{n+1-j}-1],
\end{equation*}
and
$$ \Delta U^{k,j}_{\prod_{i=1}^{n+1} p_i} \supseteq \left(\Delta V_{\prod_{i=1}^{n+1-k} p_i} \right)_{mod(\prod_{i=1}^{n+1-k}p_i)}=[0:2:\prod_{i=1}^{n-k}p_i-2  ],$$ where $j\leq \prod_{i=n-k+2}^{n+1} p_i$.

\item We can choose the integer $k$ such that all elements inside of $U^k_{\prod_{i=1}^{n+1} p_i}$ are primes. Since the range of set is large and
\begin{eqnarray*}
  \lim_{k\rightarrow \infty} \min_{k} U^{k,j}_{\prod_{i=1}^{n+1} p_i}/\{1\} &=& +\infty \\
  \lim_{k\rightarrow \infty} \max_{k} U^{k,j}_{\prod_{i=1}^{n+1} p_i}/\{1\}  &=& +\infty, \\
   &\text{and }& \\
   \Delta U^{k,j}_{\prod_{i=1}^{n+1} p_i}  &=& [0:2:\prod_{i=1}^{n-k}p_i-2  ],
\end{eqnarray*}
it leads to the infinity existence of twin prime number and to the fact that every even integer greater than 2 can be expressed as the difference of two primes

$$\liminf_{n\rightarrow \infty} \left( p_n-p_{n+1} \right)=2.$$
\item Since there is no prime between two primes that defer by 4, we can conclude the existence of infinitely many cousin primes, $p_{n+1}-p_{n}=4$, i.e.,
 $$\lim_{n\rightarrow \infty} \sum_{k>n}1_{\left\{p_k-p_{k+1}=4\right\}}\neq 0.$$
 
\item we can prove that if every even integer greater than 2 can be expressed as the difference of two primes then every even integer greater than 2 can be expressed as the sum of two primes.\\
    \item Other methods are highlighted during the proofs.
\end{itemize}

The initial idea of these above points comes from the way of generating sets of co-prime numbers recursively by adding one prime at each step as it is indicated by the following theorem

\begin{thm}  \label{theorem 1}
Any co-prime number, $n$, with the first $m$ prime numbers can be expressed according to the following equation
\begin{equation*} n= \left(\prod_{i=1}^{m}P_{i}\right) k+l,
\end{equation*}

where $k\in Z$, and $ l\in V_{mod( \prod_{i=1}^{m}P_{i})}$ a set of co-prime numbers with the first $m-1$ prime numbers, $p_2,...,p_m$.

The set of co-prime numbers, $V_{mod\left ( \prod_{i=1}^{m+1}P_{i} \right )}$, can be generated recursively as follows:

\begin{equation*} V_{mod\left ( \prod_{i=1}^{m+1}P_{i} \right )}=\left(P_{m+1}V_{mod\left ( \prod_{i=1}^{m}P_{i} \right )}+\left ( \prod_{i=1}^{m}P_{i} \right )\left [ 1:P_{m+1}-1 \right ]\right)_{mod\left ( \prod_{i=1}^{m+1}P_{i} \right )}.\end{equation*}

Concluding that any number $ x \in V_{mod( \prod_{i=1}^{m}P_{i})} \cap   \left(1:p_{m+1}^{2} \right) $ is a prime number.
\end{thm}
\
\newline

Some of important properties of set of co-prime numbers that will be very often used along this paper are highlighted by the following lemma

\begin{lem} \label{lemmaaa 1}
i. For all co-primes $x$ and $y$, we have the following two properties:
\begin{equation*}
y=\left(x\left [ 1:y-1 \right ]\right)_{mod\left ( y \right )}=\left [ 1:y-1 \right ]_{mod\left ( y \right )}
\end{equation*}
$$ \forall z_1 \in [1:y-1], \exists !  z_2 \in [1:y-1], z_2\neq z_1, \text { such that } z_1=\left(z_2 x\right)_{mod(y)}$$
ii. If $y$ is a co-prime number with $P_{i}$, $i\leq m$, we have the following properties:

$$V_{mod\left ( \prod_{i=1}^{m}P_{i} \right )}=\left(yV_{mod\left ( \prod_{i=1}^{m}P_{i}
\right )} \right)_{mod\left ( \prod_{i=1}^{m}P_{i} \right )}.$$
 $$\forall z_1 \in V_{mod\left ( \prod_{i=1}^{m}P_{i} \right )}, \exists !  z_2 \in V_{mod\left ( \prod_{i=1}^{m}P_{i} \right )}, z_2\neq z_1, \text { such that } z_1=\left(z_2 x\right)_{mod\left ( \prod_{i=1}^{m}P_{i} \right )}.$$

\end{lem}
\
\newline

The main ideas of our proofs are on creating and on analysing of set sequence of distinct co-primes with the first $n$ primes, $H=(x_1,...,x_k)$, by using Theorem \ref{theorem 1} and the important properties of the linear combination of the co-prime sets, $H=\{1,p_{n+1},..., \Pi_{i=1}^n p_i-1 \}$, that give sets of even numbers $\{0,2,4,..., \Pi_{i=1}^n p_i -2 \}$. The next three lemmas will explain our approach gradually.\\

First of all, some properties of the difference between sets of co-prime numbers will be pointed out and a recursive formula between these sets will be expressed as well.

\begin{lem} \label{lemma 1}

For all $n\geq 2$ and for all $k\geq0$, we have the following equalities

\begin{enumerate}

	\item The difference co-primes of the sets $V_.$ generate compact subsets of even numbers
	\begin{eqnarray*}
	 \left( \Delta V_{\prod_{i=1}^{n} p_i} \right)_{mod\left(\prod_{i=1}^{n} p_i\right)}&=&\left(V_{\prod_{i=1}^{n} p_i} - V_{\prod_{i=1}^{n} p_i}\right)_{mod\left(\prod_{i=1}^{n} p_i\right)}\\
	&=& \left[0:2:\prod_{i=1}^{n} p_i-2 \right].
	\end{eqnarray*}
	
	\item  $V_{\prod_{i=1}^{n+1}p_i}=\left(\left(\prod_{i=n-k+2}^{n+1}p_i \right) \left(V_{\prod_{i=1}^{n-k+1} p_i}\right) + w_k \left(\prod_{i=1}^{n-k+1}p_i \right)\right)_{mod(\prod_{i=1}^{n+1}p_i)}$,\\
 where the set $w_k$ is the set of co-primes number with $p_i$, for $ i\in \{n-k+2,...,n+1\}.$\\ \\
   The recursive relation between sets $w_.$ can de deducted, for $k>1$, as follows
   $$w_{k+1} =w_k \,p_{n-k+1}+  \left(\prod_{i=n-k+2}^{n+1}p_i\right) [1:p_{n-k+1}-1].$$
	
	\end{enumerate}

\end{lem}			
\
\newline	
As initial investigation inside the set $V_.$, a new question is raised about finding subsets that verified the first part of lemma \ref{lemma 1}, the next lemma gives interesting result

\begin{lem} \label{lemma 2}
For all $n$ bigger than 1, let's define the following sets
$$ U_{\prod_{i=1}^{n+1} p_i} = \left\{ x \in V_{\prod_{i=1}^{n+1} p_i}\; :\; x\in \left[ 0, 4 {\prod_{i=1}^{n} p_i}\right]  \right\},$$
then $$ \left(\Delta U_{\prod_{i=1}^{n+1} p_i}\right)_{mod\left(\prod_{i=1}^{n+1}p_i\right)}\supseteq \left( \Delta V_{\prod_{i=1}^{n} p_i} \right)_{mod(\prod_{i=1}^{n}p_i)},$$
and
$$ \left(\Delta U_{\prod_{i=1}^{n+1} p_i}\right) \supseteq \left[0:2:\left( \prod_{i=1}^{n}p_i\right) -2 \right] $$

\end{lem}

In order to explain the idea on how to prove that the difference or summation between all prime numbers is equals to the set of all even integers, the following lemma, a generalization and improved version of lemma \ref{lemma 2}, will help us to limit the difference only between primes numbers.

\begin{lem} \label{LemmaP}
For all integer $n>3$, it exists $k_0 \in \N$, such that for all integer $k<k_0$, the difference between co-primes numbers, $V_.$, can be expressed recursively in function of difference between subsets of $V_.$ as it is indicated below
\begin{equation*}
\left(\Delta V_{\prod_{i=1}^{n+1} p_i}\right)_{mod(\prod_{i=1}^{n+1} p_i)}=\left(\Delta U^{k,j}_{\prod_{i=1}^{n+1} p_i}\right)_{mod(\prod_{i=1}^{n+1-k} p_i)}+\sum_{j=0}^{k-1} \left(\prod_{i=1}^{n-j} p_i \right)[0:p_{n+1-j}-1],
\end{equation*}
where $U^{k,j}_{\prod_{i=1}^{n+1} p_i}=V_{\prod_{i=1}^{n+1} p_i} \bigcap  \left[j \prod_{i=1}^{n+1-k} p_i: (j+1)\prod_{i=1}^{n+1-k} p_i \right].$\\

To extend this result and omit the modulo from the right side of the previous equation i.e., from the set $\Delta U^k_{\prod_{i=1}^{n+1-k} p_i}$, we enlarge the definition of the set as $$U^{k,j}_{\prod_{i=1}^{n+1} p_i}=V_{\prod_{i=1}^{n+1} p_i} \bigcap  \left[j \prod_{i=1}^{n+1-k} p_i: (j+2) \prod_{i=1}^{n+1-k} p_i \right],$$

then it leads, under the condition that  $\Pi_{i=1}^{k-2}p_i>p_{n+2}$ to

\begin{itemize}
	\item $ \Delta U^k_{\prod_{i=1}^{n+1} p_i} \supseteq \left(\Delta V_{\prod_{i=1}^{n+1-k} p_i} \right)_{mod(\prod_{i=1}^{n+1-k}p_i)}=[0:2:\left(\prod_{i=1}^{n+1-k}p_i\right)-2  ],$
	\item For all integers $j\leq \prod_{i=n+2-k}^{n+1} p_i$, we have
$$ \Delta U^{k,j}_{\prod_{i=1}^{n+1} p_i} \supseteq \left(\Delta V_{\prod_{i=1}^{n+1-k} p_i} \right)_{mod(\prod_{i=1}^{n+1-k}p_i)}=[0:2:\left(\prod_{i=1}^{n+1-k}p_i\right)-2  ].$$
\end{itemize}

The threshold $k_0$ will be the smallest integer verified the condition, i.e.,
 $$k_0=min \left\{k:\;\Pi_{i=1}^{k-2}p_i>p_{n+2}\right\}.$$
\end{lem}

%
%
An important result will be needed in order to prove the lemma \ref{LemmaP}, evaluation the density of co-primes with $p_i, i> n$, and are less or equal to $x$, noted $\pi(x,n)$. RiChard Warlimont, \cite{coprimePSI}, and D. A. Golden, \cite{coprimePSI2}, have worked on the estimation of co-primes density under certain conditions, our following lemma will express an asymptotic density of co-primes with consecutive primes.\\

\begin{lem} \label{density_coprimes}
Let $Q$ be a set of primes and we denote by $\Psi(x,y,Q)$ the number of all positive integers less than $x$ which do not have a prime divisor which is less than $y$ and belongs to $Q$.\\
Put $P_Q:=\prod_{p_i \in Q} p_i$, $y < \min_{i}\left\{ p_i\,|\, p_i\in Q\right\}$, and let $R(y)=\frac{1}{\ln(y)}$.\\

Then we have the following statements
\begin{itemize}
	\item [(i)] For any positive integers $k$ we have $$\Psi(x,y,Q)=Card \left\{i \in (y+k P_Q, x+k P_Q] \,|\, gcd(i,P_Q)=1 \right\}.$$
	\item [(ii)] For any positive integers $k$ we have $$\Psi(x,y,Q)=\Psi(x+k P_Q,y,Q)-\Psi(y+k P_Q, y, Q)+y,$$
	and if $x<P_Q$ and $gcd(x,P_Q)=1$ then
	$$ \Psi(x,y,Q)+\Psi(P_Q-x,y,Q)=\Psi(P_Q,y,Q)+1=P_Q \prod_{p_i \in Q} \left(1-\frac{1}{p_i}\right)+1.$$
	\item [(iii)] Uniformly in $1.5\leq y\leq x$ we have
	$$\Psi(x,y,Q)=(x-y) \prod_{p_i \in Q} \left(1-\frac{1}{p_i}\right) \left( 1+O\left(R(y)\right) \right)+y.$$
	\item [(iv)] If $x$ and $y$ are big enough and $Q=\left\{p_i\,|\, y<p_i \leq x \right\}$ then
	$$\Psi(x,y,Q)\sim (x-y) \frac{\ln(y)}{\ln(x)}+y,$$
	and if $y/x\rightarrow 0$ then
	$$\Psi(x,y,Q)\sim \pi(x,y/\ln(y)) \sim x \frac{\ln(y)}{\ln(x)}.$$
\end{itemize}



\end{lem}

The coming sections will be devoted to prove the Theorem \ref{theorem 1} and the Lemma \ref{lemmaaa 1} at first step of our approach as an introduction to the proof of the lemma \ref{lemma 1}, \ref{lemma 2}, and \ref{LemmaP} respectively.\\

Finally the proof of the two conjectures with their extension will be at the end of this paper.

\section{Proof of Theorem \ref{theorem 1}}
Let $n$ and $b$ be two positive integers. If $n$ is not divisible by $a$ , then $a$ needs to be a multiple of $k$ added to another positive integer $b$ , where $b$ is smaller than $a$ ; under the condition that $a$ and $b$ are co-primes, the positive integer $n$ is co-prime with $a$. This idea leads to the generation of a large prime number from a sequence of previous prime numbers. This idea can be written in mathematical form as follows: \[n=ak+b, \text{where    } b< a, \]

if $a$ and $b$ are co--primes than $a$ and $n$ are also co-prime numbers.\\

We will generalize this idea to generate co-prime numbers recursively by including one prime at each step.\\

Let $\omega(m)$ be a set of co-primes with the first $m$ prime numbers, so $\omega$ can be expressed as: \[\omega(m)=\left ( \prod_{i=1}^{m}p_{i} \right ).k+V_{mod\left ( \prod_{i=1}^{m}p_{i} \right )},\]

A recursive question will be raised in how to find the set $V_{mod\left ( \prod_{i=1}^{m+1}p_{i} \right )}$ from the $V_{mod\left ( \prod_{i=1}^{m}p_{i} \right )}$ in order to extract all elements of $\omega(m+1)$. To do that, let suppose $\omega(m)$ is co-prime also with $\left ( p_{m+1} \right )$ and by definition it is co-prime with  $\left ( \prod_{i=1}^{m}p_{i} \right )$, so it exists two integers, $k$ and $k^{'}$ such that

\begin{equation}
\begin{cases}
\omega(m+1)=(\prod_{i=1}^{m}p_{i}).k+V_{mod( \prod_{i=1}^{m}p_{i})}\\
\omega(m+1)=( p_{m+1}).k^{'}+[ 1:( p_{m+1}-1)],
\end{cases}
\end{equation}

Let suppose the integer $x_{1}$ as a co-prime with all $ p_{i}$ , $i\leq m$, and $x_{2}$ is co-prime with $p_{m+1}$, by using Lemma  \ref{lemmaaa 1} the previous system of equations can be formulated as follows:

\begin{equation}
\begin{cases}
\omega(m+1)=(\prod_{i=1}^{m}p_{i}).k+V_{mod( \prod_{i=1}^{m}p_{i})} x_{1}\\
\omega(m+1)=( p_{m+1}).k^{'}+[ 1:( p_{m+1}-1)] x_{2},
\end{cases}
\end{equation}

Without losing generality let $x_1=p_{m+1}$ and $x_2=\prod_{i=1}^{m}p_{i}$  so $\omega(m+1)$ can be expressed by

\begin{equation}
\begin{cases}
\omega(m+1)=(\prod_{i=1}^{m}p_{i}).k+V_{mod( \prod_{i=1}^{m}p_{i})} p_{m+1}\\
\omega(m+1)=( p_{m+1}).k^{'}+[ 1:( p_{m+1}-1)] \prod_{i=1}^{m}p_{i},
\end{cases}
\end{equation}

The set of co-primes can be expressed as given below:
\begin{equation}\label{system3}
  \left ( \prod_{i=1}^{m}p_{i} \right ).k-\left ( p_{m+1} \right ).k^{'}=\left [ 1:\left ( p_{m+1} -1\right ) \right ] \prod_{i=1}^{m}p_{i}-V_{mod\left ( \prod_{i=1}^{m}p_{i} \right )} p_{m+1}
\end{equation}

Let $l_1\in \left[ 1:\left ( p_{m+1} -1\right ) \right ]$ and $l_2 \in V_{mod\left ( \prod_{i=1}^{m}p_{i} \right )}$.

To calculate $\omega(m+1)$ the equation \ref{system3} needs to be solved, by subtraction their two equations we can find
$$\prod_{i=1}^{m}p_{i}.\left ( k-l_1 \right )=\left ( p_{m+1} \right ).\left ( k^{'}-l_2 \right ),$$

where $( k-l_1)$ and $( k^{'}-l_2)$ can be calculated as:

 \[\begin{cases}
\left ( k-l_1 \right )=\vartheta .p_{m+1} \\
\left ( k^{'}-l_2 \right )=\vartheta .\left ( \prod_{i=1}^{m}p_{i} \right ),\end{cases}\]
where $\vartheta$ is any value  inside $Z$.\\

By isolating the two variables $k$ and $k^{'}$ , we have:
$$\begin{cases}
k=\vartheta .p_{m+1}+l_1. \\
k^{'}=\vartheta .\left ( \prod_{i=1}^{m}p_{i} \right )+l_2 \ ,\end{cases}$$

Substituting the value of $k$ in system (3), we have:

\[\begin{cases} \omega(m+1)=\vartheta \left ( \prod_{i=1}^{m+1}p_{i} \right )+l_1 \left ( \prod_{i=1}^{m}p_{i} \right )+V_{mod\left ( \prod_{i=1}^{m}p_{i} \right )} \\
\omega(m+1)=\vartheta \left ( \prod_{i=1}^{m+1}p_{i} \right )+l_2 \left ( p_{m+1} \right )+\left [ 1:\left ( p_{m+1}-1 \right ) \right ],\end{cases}\]

By substitution the values of $l_1$ and $l_2$ in the above equations we conclude

 \[\begin{cases}
\omega(m+1)=\vartheta .\left ( \prod_{i=1}^{m+1}p_{i} \right )+\left ( \left [ 1:\left ( p_{m+1}-1 \right ) \right ].\left ( \prod_{i=1}^{m}p_{i} \right ) \right )+\left ( p_{m+1} \right ).V_{mod\left ( \prod_{i=1}^{m}p_{i} \right )} \\

\omega(m+1)=\vartheta .\left ( \prod_{i=1}^{m+1}p_{i} \right )+\left ( \left [ 1:\left ( p_{m+1}-1 \right ) \right ].\left ( \prod_{i=1}^{m}p_{i} \right ) \right )+\left ( p_{m+1} \right ).V_{mod\left ( \prod_{i=1}^{m}p_{i} \right )},
\end{cases}\]
where both equations are equivalent due to the similarity of property of the modulo of the last term.

The general formula of $\omega(m+1)$ that are co-primes with the first $m+1$ prime numbers will be as follows:
\[\omega(m+1)=\vartheta .\left ( \prod_{i=1}^{m+1}p_{i} \right )+\left ( \left [ 1:\left ( p_{m+1}-1 \right ) \right ].\left ( \prod_{i=1}^{m}p_{i} \right ) \right )+\left ( p_{m+1} \right ).V_{mod\left ( \prod_{i=1}^{m}p_{i} \right )},\]

We can general also the set of co-primes with the first $m+1$ prime numbers in order to get the following:
\[\omega(m+1)=\vartheta .\left ( \prod_{i=1}^{m+1}p_{i} \right )+V_{mod\left ( \prod_{i=1}^{m+1}p_{i} \right )} \ ,\]

where
$$V_{mod\left ( \prod_{i=1}^{m+1}p_{i} \right )}=\left ( \left [ 1:\left ( p_{m+1}-1 \right ) \right ].\left ( \prod_{i=1}^{m}p_{i} \right ) +\left ( p_{m+1} \right ).V_{mod\left ( \prod_{i=1}^{m}p_{i} \right )} \right )_{mod\left ( \prod_{i=1}^{m+1}p_{i} \right )}.$$

\begin{flushright}
$\blacksquare$
\end{flushright}
\section{ Proof of Lemma \ref{lemmaaa 1}}
\subsection{Proof of Lemma \ref{lemmaaa 1} part (i)}
We will start proving the first part of the Lemma \ref{lemmaaa 1} with following equality:
\[\left ( \prod_{i=1}^{m}p_{i} \right ).\left [ 1:\left ( p_{m+1}-1 \right ) \right ]_{mod\left ( p_{m+1} \right )}=\left [ 1:\left ( p_{m+1}-1 \right ) \right ]_{mod\left ( p_{m+1} \right )} \ .\] \\

Since $ (\prod_{i=1}^{m}p_{i})$ is co-prime with $p_{m+1}$, there exists $k\in N$, where $(\prod_{i=1}^{m}p_{i})$ can be written as:

$( \prod_{i=1}^{m}p_{i})=k.p_{m+1}+r$,
where $1\leq r< p_{m+1}$.

Then it implies that

 $$ \left(k.p_{m+1}+r\right) \left [ 1,2,3,........\left ( p_{m+1}-1 \right ) \right ]_{p_{m+1}}=\left [ r,2r,3r,........\left ( p_{m+1}-1 \right ).r \right ]_{p_{m+1}}.$$

By contradiction we can prove that
\begin{equation}\label{equation2}
  \left [ ir \right ]_{mod\left ( p_{m+1} \right )}\neq 0
\end{equation}

If $[ ir]_{mod( p_{m+1})}= 0$, then it exists $\mu\in Z$ such that $ir=\mu p_{m+1},$ which it means that $i$ or $r$ are multiple of $p_{m+1}$, where it contradicts of the fact that $i$ and $r$ are both less than $p_{m+1}$.

To finish the proof of Lemma \ref{lemmaaa 1} part (i) and first part of (ii), we still need to show by contraction the following:
\begin{equation}\label{equation3}
  \left [ ir \right ]_{mod\left ( p_{m+1} \right )}=\left [ jr \right ]_{mod\left ( p_{m+1} \right )}, \ \ \forall \  i\neq j \in \left\{1,2..., p_{m+1}-1\right\}
\end{equation}

Let's suppose that it exists two integers $i$ and $j$ such that:
\[\begin{cases}
ir=k^{'}.p_{m+1}+r_{0},\ \  k^{'}\in \Z  \\
jr=k^{''}.p_{m+1}+r_{0},\ \  k^{''}\in \Z \ .
\end{cases}\]

By subtracting the two equation, we have that
$$( i-j)r=( k^{'}-k^{''}).p_{m+1},$$

which it contradicts with the equation \ref{equation2} due to the facts that
\begin{eqnarray*}
  r &\in& [1:p_{m+1}-1] \\
  \left(i-j\right)_{mod(p_{m+1})} &\leq&  p_{m+1}-1
\end{eqnarray*}

From the equations \ref{equation2} and \ref{equation3}, The following equation can be deducted easily
$$\left [ r \ 2r \ 3r \ ........\left ( p_{m+1}-1 \right )r \right ]_{mod\left ( p_{m+1} \right )}=\left [ 1:\left ( p_{m+1}-1 \right ) \right ]_{mod\left ( p_{m+1} \right )},$$
then it implies that

$$\left ( \prod_{i=1}^{m}p_{i} \right ).\left [ 1:\left ( p_{m+1}-1 \right ) \right ]_{mod\left ( p_{m+1} \right )}=\left [ 1:\left ( p_{m+1}-1 \right ) \right ]_{mod\left ( p_{m+1} \right )}.$$

Since for all all two co-primes,$ k_{1}$ and $k_{2}$, we can write $\left(k_1\right)_mod(k_2)=r$, where $r \in [1:k_2$, then we can Generalize our result as follows

$$\left(k_{1}\left [ 1:\left ( k_{2}-1 \right ) \right ]\right)_{mod\left ( k_{2} \right )}=\left [  1:\left ( k_{2}-1 \right )\right ]_{mod\left ( k_{2} \right )}.$$
\begin{flushright}
$\blacksquare$
\end{flushright}
\subsection{Proof of Lemma \ref{lemmaaa 1} part (ii)}
Here we will prove the second part of Lemma \ref{lemmaaa 1}
\[p_{m+1}.V_{mod\left ( \prod_{i=1}^{m}p_{i} \right )}=V_{mod\left ( \prod_{i=1}^{m}p_{i} \right )},\]

For all $m>1$; the cadinality of the set $V_.$ can be calculated as follows

$$ Card \left(V_{mod\left ( \prod_{i=1}^{m}p_{i} \right )}\right)=\prod_{i=2}^{m} (p_i-1).$$
Let $Card \left(V_{mod\left ( \prod_{i=1}^{m}p_{i} \right )}\right)=M$, so
\[V_{mod\left ( \prod_{i=1}^{m}p_{i} \right )}=\left \{ r_{1},r_{2},r_{3},......,r_{M} \right \} . \]

Now, we will test first the non existing of $i \in \{1,2...,M\}$ such that:
\begin{equation*}
r_{i}:\;\; r_{i}  p_{m+1}=k\left(\prod_{i=1}^{m}p_{i}\right),\end{equation*}
if the last equation is true, then for all $j\leq m$ we have
 $$\frac{r_{i}}{p_{j}}\in \N,$$

which it implies that
$$r_{i}=k^{'}( \prod_{i=1}^{m}p_{i}),$$

this statement contradict with the fact that $r_i \in V_{mod(\prod_{i=1}^{m}p_{i})}$, i.e.,
 $$\left( r_{i}\right)_{mod(\prod_{i=1}^{m}p_{i})}\neq 0.$$\\

To complete the proof, second test needs to be checked :

$$  \exists \text{  }i\neq j \text{ such that }\left [ \left ( r_{i} \right )p_{m+1} \right ]_{mod\left ( \prod_{i=1}^{m}p_{i} \right )}=\left [ \left ( r_{j} \right )p_{m+1} \right ]_{mod\left ( \prod_{i=1}^{m}p_{i} \right )},$$
which it implies that $r_{i}$ and $r_{j}$ have the following values: \[\begin{cases}
r_{i}p_{m+1}=k^{'}\left ( \prod_{i=1}^{m}p_{i} \right )+r_{0} \;\; k^{'} \in \Z \\
r_{j}p_{m+1}=k^{''}\left ( \prod_{i=1}^{m}p_{i} \right )+r_{0} \;\; k^{''}\in \Z\ ,
\end{cases}\]
where $r_0 \in [1:\prod_{i=1}^{m}p_{i}-1].$\\

By subtracting the above two equations, we can get
$$\left (r_{i}-r_{j}  \right )p_{m+1}=\left ( \prod_{i=1}^{m}p_{i} \right )\left ( k^{'}-k^{''} \right ),$$

which it means that $( r_{i}-r_{j})$ has $p_{1},p_{2},.....p_{m}$ as divisors, so as consequence we can write
 $$( r_{i}-r_{i})\geq( \prod_{i=1}^{m}p_{i}),$$

 again this statement contradict with the fact that the integers $r_.$ are less than $\prod_{i=1}^{m}p_{i}.$\\

We can conclude now
$$  \forall \text{  }i\neq j, \, \left [ \left ( r_{i} \right )p_{m+1} \right ]_{mod\left ( \prod_{i=1}^{m}p_{i} \right )}\neq \left [ \left ( r_{j} \right )p_{m+1} \right ]_{mod\left ( \prod_{i=1}^{m}p_{i} \right )}.$$

The  lemma \ref{lemmaaa 1} (ii), part two, can be concluded since the cardinality of the set $V_{\left(\prod_{i=1}^{m}p_{i} \right)} $ is preserved after multiplying it by a co-prime integer.\
\newline
\begin{flushright}
$\blacksquare$
\end{flushright}

\
\newline
\\

\section{Proof of Lemma \ref{lemma 1}}
\subsection{Proof of Lemma \ref{lemma 1} part (i)}
The proof by induction will be used to prove the first part. The case where $n=2$ is trivial to be checked, 1-1=0. For $n=3$ the co-prime set with the two primes 2 and 3 is defined by $V_{\prod_{i=1}^3p_i}=\lbrace 1,5\rbrace$, then
$$\Delta V_{\prod_{i=1}^3 p_i}=\lbrace 1,5\rbrace-\lbrace 1,5\rbrace=\lbrace 0,4,-4\rbrace$$
which implies that
$$\left( \Delta V_{\prod_{i=1}^3 p_i} \right) _{mod\left( \prod_{i=1}^3 p_i\right) }=\left\lbrace 0,2,4\right\rbrace  _{mod \left( \prod_{i=1}^3 p_i\right) }.$$

Let's assume that the formula is correct for $n$ and we need to prove the formula for $n+1$.
From the assumption that
$\left( \Delta V_{\prod_{i=1}^n p_i} \right)_{mod\left( \prod_{i=1}^n p_i\right)}=\left[ 0:2:\prod_{i=1}^n p_i-2  \right]$, the next formula needs to be concluded,

$$\left( \Delta V_{\prod_{i=1}^{n+1} p_i} \right) _{mod\left( \prod_{i=1}^{n+1} p_i\right) }=\left[ 0:2:\prod_{i=1}^{n+1} p_i-2  \right].$$

Note that, from Lemma \ref{lemmaaa 1}, that
$$\left( p_{n+1} V_{\prod_{i=1}^{n} p_i} \right) _{mod\left( \prod_{i=1}^{n} p_i\right) }=V_{\prod_{i=1}^{n} p_i}.$$

then
$$ \left( p_{n+1} \Delta V_{\prod_{i=1}^{n} p_i} \right) _{mod\left( \prod_{i=1}^{n} p_i\right) }=\left( \Delta V_{\prod_{i=1}^{n} p_i}\right)_{mod\left( \prod_{i=1}^{n}p_i\right) }.$$

Further, from theorem \ref{theorem 1}, we have that
$$ \Delta V_{\prod_{i=1}^{n+1} p_i} =\left( p_{n+1} V_{\prod_{i=1}^{n} p_i} +{\prod_{i=1}^{n} p_i} \left[1:p_{n+1}-1 \right]  \right) -\left( p_{n+1} V_{\prod_{i=1}^{n} p_i} +{\prod_{i=1}^{n} p_i} \left[1:p_{n+1}-1 \right]  \right)$$
this difference gives
$$\left( \Delta V_{\prod_{i=1}^{n+1} p_i}\right)_{mod \left( \prod_{i=1}^{n+1} p_i\right) } =\left( p_{n+1} \Delta V_{\prod_{i=1}^{n} p_i} +{\prod_{i=1}^{n} p_i} \left[2-p_{n+1}:p_{n+1}-2 \right]  \right)_{mod \left( \prod_{i=1}^{n+1} p_i\right) }.$$

And from lemma\ref{lemmaaa 1}, we have
$$\left(p_{n+1} \left[ 1:1:\prod_{i=3}^{n}p_i -1\right] \right)_{mod \left( \prod_{i=3}^{n} p_i\right)}=\left[ 1:1:\prod_{i=3}^{n}p_i -1\right] $$
by multiplying the previous equation by 2 and including a value zero to the both side, it leads to
$${0}\cup\left(2 p_{n+1} \left[ 1:1:\prod_{i=3}^{n}p_i -1\right] \right)_{mod \left( \prod_{i=1}^{n} p_i\right)}={0}\cup 2 \left[ 1:1:\prod_{i=3}^{n}p_i -1\right] $$
which it can rewritten as follows
$$\left(p_{n+1} \left[ 0:2:\prod_{i=1}^{n}p_i -2\right] \right)_{mod \left( \prod_{i=1}^{n} p_i\right)}=\left[ 0:2:\prod_{i=1}^{n}p_i -2\right]. $$

By combining these properties leads to
$$\left( \Delta V_{\prod_{i=1}^{n+1} p_i}\right)_{mod \left( \prod_{i=1}^{n+1} p_i\right) } =\left( p_{n+1}\left[0:2:\prod_{i=1}^{n}p_i-2 \right]  +\left({\prod_{i=1}^{n} p_i}\right) \left[2-p_{n+1}:p_{n+1}-2 \right]  \right)_{mod \left( \prod_{i=1}^{n+1} p_i\right) }$$
$$=\left(\left[0:2:\prod_{i=1}^{n}p_i-2 \right]  +\left({\prod_{i=1}^{n} p_i}\right) \left[2-p_{n+1}:p_{n+1}-2 \right]  \right)_{mod \left( \prod_{i=1}^{n+1} p_i\right) }.$$

The key idea of this lemma comes from previous equation, which it can be derived differently as the following\\

For all $i \in \left[ 0:2:\prod_{i=1}^n p_i-2\right], \;\; \exists l_0 <p_{n+1}$ such that
$$i+l_0 \prod_{i=1}^n p_i \in p_{n+1}\left[ 0:2:\prod_{i=1}^n p_i-2\right], $$
this yields

$$\left( \Delta V_{\prod_{i=1}^{n+1} p_i}\right)_{mod \left( \prod_{i=1}^{n+1} p_i\right) } =\left( \bigcup_{i\in \left[ 0:2:\prod_{i=1}^n p_i-2\right] } \left\lbrace  i+l_0 \prod_{i=1}^n p_i  \right\rbrace  + (\prod_{i=1}^{n} p_i) \left[2-p_{n+1}:p_{n+1}-2 \right]  \right)_{mod \left( \prod_{i=1}^{n+1} p_i\right) }. $$

Since
$$\left(    l_0 \prod_{i=1}^n p_i    + (\prod_{i=1}^{n} p_i) \left[2-p_{n+1}:p_{n+1}-2 \right]  \right)_{mod \left( \prod_{i=1}^{n+1} p_i\right) }=\left( (\prod_{i=1}^{n} p_i) \left[2-p_{n+1}:p_{n+1}-2 \right]  \right)_{mod \left( \prod_{i=1}^{n+1} p_i\right) }$$
which it leads to

$$ \left( \Delta V_{\prod_{i=1}^{n+1} p_i}\right)_{mod \left( \prod_{i=1}^{n+1} p_i\right) } =\left( \cup_{i\in \left[ 0:2:\prod_{i=1}^n p_i-2\right] } \left\lbrace  i  \right\rbrace  + (\prod_{i=1}^{n} p_i) \left[2-p_{n+1}:p_{n+1}-2 \right]  \right)_{mod \left( \prod_{i=1}^{n+1} p_i\right) }   $$
and
$$ \left( \Delta V_{\prod_{i=1}^{n+1} p_i}\right)_{mod \left( \prod_{i=1}^{n+1} p_i\right) } =\left( \left[ 0:2:\prod_{i=1}^n p_i-2\right]  + (\prod_{i=1}^{n} p_i) \left[2-p_{n+1}:p_{n+1}-2 \right]  \right)_{mod \left( \prod_{i=1}^{n+1} p_i\right) }   .$$

Finally, from the last equation, it follows that

$$ \left( \Delta V_{\prod_{i=1}^{n+1} p_i}\right)_{mod \left( \prod_{i=1}^{n+1} p_i\right) } =\left[ 0:2:\prod_{i=1}^{n+1} p_i-2\right].$$
\begin{flushright}
$\blacksquare$
\end{flushright}
\subsection{Proof of Lemma \ref{lemma 1} part (ii)}
From the recursive formula of Theorem \ref{theorem 1}, we have the following

\begin{eqnarray*}
  V_{\prod_{i=1}^{n+1}p_i}&=&\left(p_{n+1} V_{\prod_{i=1}^{n}p_i} + \prod_{i=1}^{n}p_i [1:p_{n+1}-1]\right)_{mod(\prod_{i=1}^{n+1}p_i)}\\
&=&p_{n+1} \left( p_{n} V_{\prod_{i=1}^{n-1}p_i} + \prod_{i=1}^{n-1}p_i [1:p_{n}-1]\right)+ \left(\prod_{i=1}^{n}p_i \right) [1:p_{n+1}-1]\\
&=&p_{n+1} p_{n} V_{\prod_{i=1}^{n-1}p_i}+ p_{n+1} \prod_{i=1}^{n-1}p_i [1:p_{n}-1] +\prod_{i=1}^{n}p_i [1:p_{n+1}-1]\\
&=&p_{n+1} p_{n} V_{\prod_{i=1}^{n-1}p_i}+ \prod_{i=1}^{n-1}p_i\left( p_{n+1} [1:p_{n}-1] +p_{n} [1:p_{n+1}-1]\right)
\end{eqnarray*}
Further, the set $V_{\prod_{i=1}^{n+1}p_i}$ can be written in term of $V_{\prod_{i=1}^{n-1}p_i}$ as follows

$$V_{\prod_{i=1}^{n+1}p_i}=\left(p_{n+1} p_{n} V_{\prod_{i=1}^{n-1}p_i}+ \left(\prod_{i=1}^{n-1}p_i\right)  \left( w_2\right)\right)_{mod(\prod_{i=1}^{n+1}p_i)},$$

where $w_2=\left( p_{n+1} [1:p_{n}-1] +p_{n} [1:p_{n+1}-1]\right).$
We can go further by generating general sequence in terms of $w_i$ as follows

$$V_{\prod_{i=1}^{n+1}p_i}=\left(\prod_{i=n-k+2}^{n+1}p_i \right) \left(V_{\prod_{i=1}^{n-k+1} p_i}\right) + w_k \left(\prod_{i=1}^{n-k+1}p_i \right),$$
the last equation leads to
$$V_{\prod_{i=1}^{n+1}p_i}=\prod_{i=n-k+2}^{n+1}p_i \left(p_{n-k+1} \left(V_{\prod_{i=1}^{n-k} p_i } \right)+\left(\prod_{i=1}^{n-k} p_i\right) [1:p_{n-k+1}-1] \right) + w_k \prod_{i=1}^{n-k+1}p_i ,$$
then the general term of the sequence, $w_k$, can be obtained from $w_k$ as recursive formula
$$w_{k+1} \prod_{i=1}^{n-k}p_i=w_k \prod_{i=1}^{n-k+1}p_i + \left(\prod_{i=1}^{n-k}p_i\right) \left(\prod_{i=n-k+2}^{n+1}p_i\right) [1:p_{n-k+1}-1], \text{ } \forall k\geq 1,$$
Then
$$w_{k+1} =w_k p_{n-k+1}+  \left(\prod_{i=n-k+2}^{n+1}p_i\right) [1:p_{n-k+1}-1].$$

It still remains to prove that $w_k$ is is the set of co-primes with $p_i$, where $ i=\{n-k+2,,...,n+1\}$. Let's start by $w_2$ since $w_1$ is obvious.\\
By definition $V_{p_{n+1} p_n}$ is a set of all co-prime number with $p_{n+1}$ and $p_{n}$ and bounded by $p_{n+1} p_n$ and also $V_{p_{n}}$ and $V_{p_{n+1}}$ be co-primes number with $p_{n}$ and $p_{n+1}$ respectively, then
\begin{equation}
\begin{cases}
V_{p_{n}}= kp_n +[1:p_n-1], \;\; k\in \Z \\
V_{p_{n+1}}= k'p_{n+1} +[1:p_{n+1}-1]\;\; k'\in \Z,\\
\end{cases}
\end{equation}
and From Lemma \ref{lemmaaa 1}, in mod $p_n$ equation, we have
\begin{equation}
\begin{cases}
V_{p_{n}}= kp_n +p_{n+1}[1:p_n-1]\\
V_{p_{n+1}}= k'p_{n+1} +p_n[1:p_{n+1}-1]\\
\end{cases}
\end{equation}

To get the form of co-primes with $p_{n}$ and $p_{n+1}$, $V_{p_{n}}\bigcap V_{p_{n+1}}$, then by Subtracting both equations we get

$$k'p_{n+1}-k p_n=p_{n+1}[1:p_n-1]-p_n[1:p_{n+1}-1],$$

Since $p_{n+1}[1:p_n-1]-p_n[1:p_{n+1}-1]$ is co-primes with $p_{n}$ and $p_{n+1}$, then it exists particular solution $k_0$ and $k_0^{'}$, that leads to general solution
\begin{equation}
\begin{cases}
k=k_0+h p_{n+1}\\
k'=k_0^{'}+ hp_n,
\end{cases}
\end{equation}
where $h$ is any integer.

Hence, $V_{p_n}$ can be written under the following form
$$V_{p_n}\bigcap V_{p_{n+1}}=p_{n+1} p_{n} h+p_{n} k_0 +[1:p_n-1],$$
where $k_0=k_0^{'} p_{n} + p_{n+1}[1:p_n-1]-p_n[1:p_{n+1}-1]$, then

$$V_{p_n}\bigcap V_{p_{n+1}}=p_{n+1} p_{n} (h+k_0^{'})+p_n[1:p_{n+1}-1]-p_{n+1}[1:p_n-1].$$

Since $p_{n+1}$ and $p_{n}$ are co-primes and by using the property of mod from lemma 1, then $\left(p_{n+1}p_{n}-p_{n+1}[1:p_n-1]\right)_{mod(p_{n+1}p_{n})}=\left(p_{n+1}[1:p_n-1]\right)_{mod(p_{n+1}p_{n})}$,
The set $V_{p_{n+1}p_{n}}$, co-prime numbers with $p_{n+1}p_{n}$, can be defined as follows:
\begin{eqnarray*}
V_{p_{n+1}p_{n}} &=& V_{p_n}\bigcap V_{p_{n+1}} \\
                 &=& \left(p_n[1:p_{n+1}-1]-p_{n+1}[1:p_n-1]\right)_{mod(p_{n+1}p_{n})} \\
                 &=& \left(p_n[1:p_{n+1}-1]+p_{n+1}[1:p_n-1]\right)_{mod(p_{n+1}p_{n})},
\end{eqnarray*}
then $w_k=V_{p_{n+1}p_{n}}.$\\

We can prove easily that $w_{k+1}$ is co-primes with $p_i$, where $ i=\{n-k+1,n-k+2,,...,n+1\}$ under the assumption that $w_k$ is co-primes with $p_i$, where $ i=\{n-k+2,,...,n+1\}$ by using the general formula of $w_.$

$$w_{k+1} =w_k p_{n-k+1}+  \left(\prod_{i=n-k+2}^{n+1}p_i\right) [1:p_{n-k+1}-1].$$
\begin{flushright}
  $\blacksquare$
\end{flushright}

\
\newline
\\
\section{Proof of Lemma \ref{lemma 2}}
The lemma will be proven by using induction method. Let check the formula for $n=2$

The set of co-prime numbers with the two primes 2 and 3 is defined as $V_{mod\left(2*3\right)}=\{1,5\}$ then $ \Delta V_{mod\left(2*3\right)}=\left\{0,2,4\right\}_{mod(2*3)}$.\\

From Theorem \ref{theorem 1}, we have that
$$V_{mod\left(2*3*5\right)}=\left(5*V_{mod\left(2*3\right)}+2*3 [1:4]\right)_{mod\left(2*3*5\right)}=\left\{1,7,11,13,17,19,23,29\right\}$$
Therefore
$$U_{\left\{2*3*5\right\}}=\left\{1,7,11,13,17,19,23\right\} \Rightarrow \left(\Delta U_{\left\{2*3*5\right\}}\right)_{mod\left(2*3*5\right)}=\left[0:2:6\right].$$

From lemma \ref{lemmaaa 1}, the two important equations can be extracted
\begin{equation}
\begin{cases}
V_{\prod_{i=1}^{n+1}p_i}=\left(p_{n+1}V_{\prod_{i=1}^{n}p_i}+\prod_{i=1}^{n}p_i\left[1:p_{n+1}-1\right]\right)_{mod\left(\prod_{i=1}^{n+1}p_i\right)}\\
\text{and}\\
\left(p_{n+1}V_{\prod_{i=1}^{n}p_i}\right)_{mod\left(\prod_{i=1}^{n}p_i\right)}=V_{\prod_{i=1}^{n}p_i}.
\end{cases}
\end{equation}

The set $U_{\prod_{i=1}^{n+1}p_i}$ is defined as an intersection of the set $V_{\prod_{i=1}^{n+1}p_i}$ and the interval
$\left[0, 4\prod_{i=1}^{n}p_i\right]$, i.e.,
$$U_{\prod_{i=1}^{n+1}p_i}=V_{\prod_{i=1}^{n+1}p_i} \bigcap \left[0: 4\prod_{i=1}^{n}p_i\right].$$

For all element $i$ in the set $V_{\prod_{i=1}^{n}p_i}$, it exists unique element $j$, different from  $i$, in $V_{\prod_{i=1}^{n}p_i}$ and an integer $k_0 \in \N$ such that
$$ip_{n+1}=j+k_0\prod_{i=1}^{n}p_i,$$
which it leads to the following equation
$$p_{n+1}V_{\prod_{i=1}^{n}p_i}=\bigcup_{i=0}^{p_{n+1}-1} \Omega_i,$$

where $\Omega_i=\left\{y :\; i \prod_{i=1}^{n}p_i<y<(i+1) \prod_{i=1}^{n}p_i,\; \exists j \in V_{\prod_{i=1}^{n}p_i} \text{ where } jp_{n+1}=y \right\}.$ or just simply by
$$\Omega_i=p_{n+1} V_{\prod_{i=1}^{n}p_i} \bigcap [i \prod_{i=1}^{n}p_i, (i+1)\prod_{i=1}^{n}p_i ].$$

Then the set $V_{\prod_{i=1}^{n+1}p_i}$ can be decomposed in terms of $\Omega_.$ as follows:

\begin{align}
	V_{\prod_{i=1}^{n+1}p_i}=\left(p_{n+1} V_{\prod_{i=1}^{n}p_i}+\prod_{i=1}^{n}p_i\left[1:p_{n+1}-1\right]\right)_{mod\left(\prod_{i=1}^{n+1}p_i\right)}\\
=\left(\bigcup_{i=0}^{p_{n+1}-1} \Omega_i+\prod_{i=1}^{n}p_i\left[1:p_{n+1}-1\right]\right)_{mod\left(\prod_{i=1}^{n+1}p_i\right)}\\
=\bigcup_{i=0}^{p_{n+1}-1} \left(i \left(\prod_{i=1}^{n}p_i\right)+\bigcup_{l=0, l\neq i}^{p_{n+1}-1}  \left(\Omega_l\right)_{mod\left( \prod_{i=1}^{n}p_i  \right)}  \right)_{mod\left( \prod_{i=1}^{n+1}p_i  \right)}
\end{align}

To make the notation easier to deal with, let define the following sets
$$W_i= \left(i \left(\prod_{i=1}^{n}p_i\right)+\bigcup_{l=0, l\neq i}^{p_{n+1}-1}  \left(\Omega_l\right)_{mod\left( \prod_{i=1}^{n}p_i  \right)}   \right),$$
then
$$V_{\prod_{i=1}^{n+1}p_i}=\bigcup_{i=0}^{p_{n+1}-1} W_i. $$

From the properties of the co-primes that for all element $x$ in $V_{\prod_{i=1}^{n}p_i}$ then the element $\prod_{i=1}^{n}p_i -x$ is also in $V_{\prod_{i=1}^{n}p_i}$, the below equation can be concluded

\begin{equation}\label{star1}
  W_0 \cup \left(\Omega_{0}\right)_{mod\left( \prod_{i=1}^{n}p_i  \right)}=V_{\prod_{i=1}^{n}p_i}
\end{equation}

and
$$\left(\Omega_{0}\right)_{mod\left( \prod_{i=1}^{n}p_i  \right)} = \left(\prod_{i=1}^{n}p_i- \Omega_{p_{n+1}-1} \right)_{mod\left( \prod_{i=1}^{n}p_i  \right)}$$

In general we can write

$$\left(W_i \bigcup \left(W_{i+1}-\prod_{i=1}^{n}p_i\right)\right)_{mod\left( \prod_{i=1}^{n}p_i \right)}=V_{\prod_{i=1}^{n}p_i},$$
and for any two different integers $i$ and $j$ we can write 
$$\left(W_i \bigcup \left(W_{j}-(j-i)\prod_{i=1}^{n}p_i\right)\right)_{mod\left( \prod_{i=1}^{n}p_i \right)}=V_{\prod_{i=1}^{n}p_i}.$$
Then we can write that

$$\Delta \left(V_{\prod_{i=1}^{n+1}p_i} \bigcap [0,2 \prod_{i=1}^{n}p_i] \right)\supseteq \left(\Delta \left(V_{\prod_{i=1}^{n}p_i}\right)\right)_{mod(\prod_{i=1}^{n}p_i)}.$$

In order to omit the $mod\left( \prod_{i=1}^{n}p_i  \right)$ from our previous question, it is important to deal with the negative values in different ways.

For that, we need to prove that for all $x$ and $y$ are in $V_{\prod_{i=1}^{n}p_i}$, than $\left|x-y\right|$ and $-\left|y-x\right|+\prod_{i=1}^{n}p_i $ are both in  $\left( \Delta U_{\prod_{i=1}^{n+1} p_i} \right) $.

It leads from equation \ref{star1} that $x$ and $y$ belong to at least three sets from $W_0, W_1, W_2,$ and $ W_3$, i.e.,

\begin{eqnarray*}
  \text{If } x \notin W_1 &\Rightarrow& x-\prod_{i=1}^{n}p_i \in W_0, \prod_{i=1}^{n}p_i -x \in W_1,\; x+\prod_{i=1}^{n}p_i \in W_2, \text{ and } x+2\prod_{i=1}^{n}p_i \in W_3 \\
  \text{ or } &&\\
  \text{If } y \notin W_2  & \Rightarrow & y-2\prod_{i=1}^{n}p_i \in W_0, x-\prod_{i=1}^{n}p_i \in W_1,\; \prod_{i=1}^{n}p_i -y \in W_2 \text{ and } x+\prod_{i=1}^{n}p_i \in W_3
\end{eqnarray*}
In general for $j =\{i-2,i-1, i+1, i+2\}$ we have the following implication 
$$ \text{ if }\; x +i\prod_{i=1}^{n}p_i  \notin W_i \Rightarrow x+j\prod_{i=1}^{n}p_i \in W_j, $$\\

In order to cover all cases, it is enough to consider the set $U$ as follows

$$U_{\prod_{i=1}^{n+1}p_i}=W_0 \cup W_1 \cup W_2 \cup W_3  \cup W_4$$

Now it is clear to observe that
\begin{align}
\left(\Delta U_{\prod_{i=1}^{n+1}p_i} \right)_{mod\left( \prod_{i=1}^{n+1}p_i \right)} \supseteq \left(\Delta V_{\prod_{i=1}^{n}p_i} \right)_{mod\left( \prod_{i=1}^{n}p_i \right)}\\
\left(\Delta U_{\prod_{i=1}^{n+1}p_i} \right) \supseteq \left[0:2: \prod_{i=1}^{n}p_i -2 \right]
\end{align}

\begin{flushright}
$\blacksquare$
\end{flushright}\
\newline
\\

\section{Proof of Lemma \ref{LemmaP}}
%
%
%

Let start with a small examples to make all properties clear as it is indicated by follows\\

From Theorem \ref{theorem 1}, the co-prime sets can be calculated as following
\begin{itemize}
  \item[] $V_{2*3} = \{1,5\}$,
  \item[] $V_{2*3*5} = \{1,7,11,13,17,19,23,29\}$,
  \item[] $V_{2*3*5*7} = \{1,11,13,17,19,23,29,31,37,41,43,47,59,61,67,...,197,199,209\}.$
\end{itemize}
By simple calculation we can find the following results
\begin{itemize}
  \item[] $\left(\Delta V_{2*3*5}\right)_{mod(2*3*5)}=\{0, 2, 4, 6, 8, 10, 12, 14, 16, 18, 20, 22, 24, 26, 28\}$.
  \item[] $\left(\Delta V_{2*3*5*7}\right)_{mod(2*3*5*7)}=\{0, 2, 4,..., 206, 208\}$.
  \item[] $\Delta U^1_{2*3*5*7}=\{0, 2, 4, 6, 8, 10, 12, 14, 16, 18, 20, 22, 24, 26, 28\}.$
  \item[] $\Delta U^1_{2*3*5*7*11}=\{0, 2, 4,..., 206, 208\}.$
  \item[] $\Delta U^2_{2*3*5*7*11}=\{0, 2, 4, 6, 8, 10, 12, 14, 16, 18, 20, 22, 24, 26, 28\}.$
  \item[] $\Delta U^3_{2*3*5*7*11*13}=\{0, 2, 4, 6, 8, 10, 12, 14, 16, 18, 20, 22, 24, 26, 28\}.$
  \item[] $\Delta U^5_{2*3*5*7*11*13*17*19*23}=\{0, 2, 4,..., 206, 208\}.$
\end{itemize}

To make this result general, we will proceed to prove the first part of this lemma by using the Theorem \ref{theorem 1}, for which it gives that for all $n>2$:

$$  V_{\prod_{i=1}^{n+1} p_i} = \left(p_{n+1} V_{\prod_{i=1}^{n} p_i}+\left( \prod_{i=1}^{n+1} p_i \right) \left[1:p_{n+1}-1  \right]\right)_{mod(\prod_{i=1}^{n+1} p_i)}. $$

In order to track the evolution of the co-prime sets, we will define
$$\Psi_i=V_{\prod_{i=1}^{n} p_i}\bigcap \left[i \frac{\prod_{i=1}^{n} p_i}{p_{n+1}}: (i+1) \frac{\prod_{i=1}^{n} p_i}{p_{n+1}}\right],$$

then we can write
\begin{equation} \label{equation17}
	V_{\prod_{i=1}^{n+1} p_i}= \left(p_{n+1} \left(\bigcup_{i=0}^{p_{n+1}-1} \Psi_i \right)+\left( \prod_{i=1}^{n+1} p_i \right) \left[1:p_{n+1}-1\right]\right)_{mod(\prod_{i=1}^{n+1} p_i)}
\end{equation}

Let define for instant the set
$$U^1_{\prod_{i=1}^{n+1} p_i}=V_{\prod_{i=1}^{n+1} p_i} \bigcap [0: \prod_{i=1}^{n} p_i]$$

and we will explain later why we need to enlarge the set $U^1_.$ by double.\\

We can find out that $\left(\Delta U^1_{\prod_{i=1}^{n+1} p_i}\right)_{mod(\prod_{i=1}^{n} p_i)}=\left(\Delta V_{{\prod_{i=1}^{n} p_i}}\right)_{mod(\prod_{i=1}^{n} p_i)} $ because of the following properties\\

\begin{itemize}
  \item The set $ U^1_{\prod_{i=1}^{n+1} p_i}$ contains all elements inside $V_{{\prod_{i=1}^{n} p_i}}$ except multiples of $p_{n+1}$ then $U^1_{\prod_{i=1}^{n+1} p_i} \subseteq V_{{\prod_{i=1}^{n+1} p_i}} $, and the missing numbers are only multiples of $p_{n+1}$ represented by the set $p_{n+1} \Psi_0$.
  
	\item If $x$ is multiple of $p_{n+1}$ then $ \left(\prod_{i=1}^{n} p_i\right) -x$ is not a multiple of $p_{n+1}$, so we can find out that if the value $x$ is not in $ U^1_{\prod_{i=1}^{n+1} p_i}$ and it in $V_{{\prod_{i=1}^{n} p_i}}$ then $$ \left(\prod_{i=1}^{n} p_i\right) -x \in U^1_{\prod_{i=1}^{n+1} p_i}.$$
      This result comes from the impossibility to have both $ x $ and $\left(\prod_{i=1}^{n} p_i\right) -x$ as multiple of $p_{n+1}$ which it leads to a contradiction
      $$\frac{\left(\prod_{i=1}^{n} p_i\right)}{p_{n+1}}\in \Z.$$
  
	\item From equation \ref{equation17}, we can conclude that for all $x\in V_{\prod_{i=1}^{n-1} p_i}$, it exists a unique integer $l<p_n$ such that
  $$ \left(x\right)_{mod(\prod_{i=1}^{n-1} p_i)} \notin \left( V_{\prod_{i=1}^{n} p_i}\bigcap \left[ l {\prod_{i=1}^{n-1} p_i}, (l+1) {\prod_{i=1}^{n-1} p_i} \right]\right)_{mod(\prod_{i=1}^{n-1} p_i)}.$$
  This result can be found out from the cardinality of the sets 
	$$\text{Card} (V_{\prod_{i=1}^{n-1} p_i})=\text{Card} (V_{\prod_{i=1}^{n-2} p_i}) (p_{n-1}-1)=\prod_{i=2}^{n} (p_i-1),$$
	or from the theorem \ref{theorem 1} and Lemma \ref{lemmaaa 1}, i.e.,
  $$\left\{x+l \prod_{i=1}^{n-1} p_i\,:\; l=0,1,...,p_n-1\right\}$$
  has at most one element which it is a multiple of $ p_n$.
\end{itemize}\
\newline

By using the above properties, we can prove that for all $x$ and $y$ inside the set $V_{{\prod_{i=1}^{n} p_i}}$ then we can prove easily that $(x-y)_{mod(\prod_{i=1}^{n} p_i)}$ and $(y-x)_{mod(\prod_{i=1}^{n} p_i)}$ are both inside $\left(\Delta U^1_{\prod_{i=1}^{n+1} p_i}\right)_{mod(\prod_{i=1}^{n})}$. To do that it is enough to consider only the cases where $x$ or $y$ are multiple of $p_{n+1}$, i. e. belong to $\Psi_0 p_{n+1}$. All cases are checked as it is indicated below
\begin{description}
  \item[Case 1] If $x \notin U^1_{\prod_{i=1}^{n+1} p_i} $, $y \in U^1_{\prod_{i=1}^{n+1} p_i}$, and $\prod_{i=1}^{n} p_i-y \in U^1_{\prod_{i=1}^{n+1} p_i} $ then $\prod_{i=1}^{n} p_i-x \in U^1_{\prod_{i=1}^{n+1} p_i}$ so we have
 $$\left((\prod_{i=1}^{n} p_i-x)-(\prod_{i=1}^{n} p_i-y)\right)_{mod(\prod_{i=1}^{n} p_i)}=(y-x)_{mod(\prod_{i=1}^{n} p_i)}$$
 $$\left((\prod_{i=1}^{n} p_i-y)-(\prod_{i=1}^{n} p_i-x)\right)_{mod(\prod_{i=1}^{n} p_i)}=(x-y)_{mod(\prod_{i=1}^{n} p_i)}.$$

 \item[Case 2] If $x \notin U^1_{\prod_{i=1}^{n+1} p_i} $, $y \in U^1_{\prod_{i=1}^{n+1} p_i}$, and $\prod_{i=1}^{n} p_i-y \notin U^1_{\prod_{i=1}^{n+1} p_i} $ then from decomposition of $V_{\prod_{i=1}^{n}p_i}$ in function of $V_{\prod_{i=1}^{{n-1}}p_i}$ so it exist an integer $l<p_{n}$ such that $x+l \prod_{i=1}^{n-1} p_i$ and $y+l \prod_{i=1}^{n-1} p_i$, in modulo $\prod_{i=1}^{n} p_i$, are inside $U^1_{\prod_{i=1}^{n+1} p_i} $ then
 $$\left((l \prod_{i=1}^{n-1} p_i+y)-(l \prod_{i=1}^{n-1} p_i+x)\right)_{mod(\prod_{i=1}^{n} p_i)}=(y-x)_{mod(\prod_{i=1}^{n} p_i)}$$
 $$\left((l \prod_{i=1}^{n-1} p_i+x)-(l \prod_{i=1}^{n-1} p_i+y)\right)_{mod(\prod_{i=1}^{n} p_i)}=(x-y)_{mod(\prod_{i=1}^{n} p_i)}.$$

This result can be justified by the fact that it not possible to find two integers less than $p_{n}$, $l_1$ and $l_2$, such that
$$\left( l_1 \prod_{i=1}^{n-1} p_i +x \right)_{mod(p_j)}=\left( l_2 \prod_{i=1}^{n-1} p_i +x  \right)_{mod(p_{j})}=0,$$
where $j$ equals to $n$ or $n+1$.

  \item[Case 3] If $x \notin U^1_{\prod_{i=1}^{n+1} p_i} $ and $y \notin U^1_{\prod_{i=1}^{n+1} p_i}$ then $\prod_{i=1}^{n} p_i-x \in U^1_{\prod_{i=1}^{n+1} p_i} $ and $\prod_{i=1}^{n} p_i-y \in U^1_{\prod_{i=1}^{n+1} p_i}$ so again
 $$\left((\prod_{i=1}^{n} p_i-x)-(\prod_{i=1}^{n} p_i-y)\right)_{mod(\prod_{i=1}^{n} p_i)}=(y-x)_{mod(\prod_{i=1}^{n} p_i)}$$
 $$\left((\prod_{i=1}^{n} p_i-y)-(\prod_{i=1}^{n} p_i-x)\right)_{mod(\prod_{i=1}^{n} p_i)}=(x-y)_{mod(\prod_{i=1}^{n} p_i)}.$$
\end{description}\
\newline

After what it said above, we can conclude that
\begin{eqnarray*}
  \left(\Delta U^1_{\prod_{i=1}^{n+1} p_i}\right)_{mod(\prod_{i=1}^{n} p_i)} &=& \left(\Delta V_{\prod_{i=1}^{n} p_i}\right)_{mod(\prod_{i=1}^{n} p_i)} \\
   &=& [0:2:\prod_{i=1}^{n}p_i  -2].
\end{eqnarray*}

By the same way we can find also that for all $j <p_{n+1}$
$$\left(\Delta U^{1,j}_{\prod_{i=1}^{n+1} p_i}\right)_{mod(\prod_{i=1}^{n} p_i)}=\left(\Delta V_{\prod_{i=1}^{n} p_i}\right)_{mod(\prod_{i=1}^{n} p_i)},$$
where $U^{1,j}_{\prod_{i=1}^{n+1} p_i}=V_{\prod_{i=1}^{n+1} p_i} \bigcap [j \prod_{i=1}^{n} p_i: (j+1)\prod_{i=1}^{n} p_i].$ \\

By using the same approach, the first part of our lemma can be derived when $x \in V_{\prod_{i=1}^{n} p_i}$ and $y \in l\prod_{i=1}^{n} p_i+ V_{\prod_{i=1}^{n} p_i}$, where $l$ is any integer number less than $p_{n+1}$
$$\left(\Delta V_{\prod_{i=1}^{n+1} p_i}\right)_{mod(\prod_{i=1}^{n+1} p_i)}=\left(\Delta U^1_{\prod_{i=1}^{n+1} p_i}\right)_{mod(\prod_{i=1}^{n} p_i)}+\left(\prod_{i=1}^{n} p_i \right)[0:1:p_{n+1}-1].$$

We suggest to prove the formula for $k=2$ in order to make the generalization easier.  Let's start by highlighting some extremely important properties, which they are derived from the symmetry of the set $V_{\prod_{i=1}^{n-1} p_i}$, i.e.,

$$\forall x\in V_{\prod_{i=1}^{n-1}p_i}\Rightarrow \prod_{i=1}^{n-1}p_i-x\in V_{\prod_{i=1}^{n-1}p_i}$$
and from the recurrence equation $$V_{\prod_{i=1}^{n} p_i} = p_{n} V_{\prod_{i=1}^{n-1} p_i}+\left( \prod_{i=1}^{n-1} p_i \right) \left[1:p_{n}-1  \right].$$\\

These properties link between elements of $V_{\prod_{i=1}^{n-1} p_i}$ and $U^2_{\prod_{i=1}^{n+1} p_i}$ as follows

For all $x \in V_{\prod_{i=1}^{n-1}p_i}$, it exists an integer set $L$ such that for all $l \in L$ we have

$$\left(l \left(\prod_{i=1}^{n-2} p_i\right) + x_{mod(\prod_{i=1}^{n-2}p_i)}\right)_{mod(\prod_{i=1}^{n-1}p_i)} \in U^2_{\prod_{i=1}^{n+1} p_i},$$

where $\text Card (L \cap [0:1:p_{n-1}-1]\geq p_{n-1}-2.$\\

This statement comes from the fact that, for all integer $i\in\{0,1,2\}$,  it is not possible to find more than one integer from the set
$$\Upsilon^2_n=\left\{l \left(\prod_{i=1}^{n-2} p_i\right) + x : l<p_{n-1}\right\},$$
which are multiple of $p_{n-1+i}$. i.e.,

\begin{equation}\label{property2}
  \text{card}\left(\Upsilon^2_n \bigcup \{I\in \N : I_{mod(p_{n-1+i})}=0\}\right)\leq 1.
\end{equation}
\\

Therefore, by the pigeonhole principle we can conclude the following
$$ \forall (x,y) \in \left(V_{\prod_{i=1}^{n-1} p_i}\right)^2 \Longrightarrow \exists \, l \in \Z:\, \left\{ x+l\prod_{i=1}^{n-2} p_i,\,y+l\prod_{i=1}^{n-2} p_i \right\}_{mod(\prod_{i=1}^{n-1} p_i)}\subseteq U^2_{\prod_{i=1}^{n+1} p_i}, $$

which implies that
\begin{eqnarray*}
  \left(\Delta U^2_{\prod_{i=1}^{n+1} p_i}\right)_{mod(\prod_{i=1}^{n-1} p_i)} &=& \left(\Delta V_{\prod_{i=1}^{n-1} p_i}\right)_{mod(\prod_{i=1}^{n-1} p_i)} \\
   &=& \left(\Delta U^{2, j}_{\prod_{i=1}^{n+1} p_i}\right)_{mod(\prod_{i=1}^{n-1} p_i)} \\
   &=& [0:2:\prod_{i=1}^{n-1} p_i-2]
\end{eqnarray*}
and also for all $j\in \{0,1,...,p_np_{n_1}-1\}$
$$\left(\Delta U^2_{\prod_{i=1}^{n+1} p_i}\right)_{mod(\prod_{i=1}^{n-1} p_i)}= \left(\Delta U^{2,j}_{\prod_{i=1}^{n+1} p_i}\right)_{mod(\prod_{i=1}^{n-1} p_i)}.$$
Then we can write
\begin{eqnarray*}
  \left(\Delta V_{\prod_{i=1}^{n+1} p_i}\right)_{mod(\prod_{i=1}^{n+1} p_i)} &=& \left(\Delta U^2_{\prod_{i=1}^{n+1} p_i}\right)_{mod(\prod_{i=1}^{n-1} p_i)}+\left(\prod_{i=1}^{n-1} p_i \right)[0:1:p_{n}-1] \\
   &+& \left(\prod_{i=1}^{n} p_i \right)[0:1:p_{n+1}-1].
\end{eqnarray*}\
\newline

By the same way we can generalize our above result to conclude that
\begin{eqnarray*}
  \left(\Delta U^{k+1}_{\prod_{i=1}^{n+1} p_i}\right)_{mod(\prod_{i=1}^{n-k} p_i)} &=& \left(\Delta V_{\prod_{i=1}^{n-k} p_i}\right)_{mod(\prod_{i=1}^{n-k} p_i)} \\
   &=& \left(\Delta U^{k+1, j}_{\prod_{i=1}^{n+1} p_i}\right)_{mod(\prod_{i=1}^{n-k} p_i)} \\
   &=& [0:2:\prod_{i=1}^{n-k} p_i-2],
\end{eqnarray*}

 where the value of $k$ should verify the following condition
\begin{equation}\label{condition1}
  k+1<\frac{p_{n-k}}{2}.
\end{equation}

This condition comes from the facts that
\begin{itemize}
  \item For all $j \in \{0,1,2..,k, k+1\}$ and for all $x \in \prod_{i=1}^{n-k} p_i$, we have the following fact
	$$\text{Card}\left(\Upsilon^{k+1}_n \bigcap \{J \in \N: J_{mod(p_{n-k+j})}=0\}\right)\leq 1,$$
	where $x_0=x_{mod(\prod_{i=1}^{n-k-1} p_i)}$ and $\Upsilon^{k+1}_n=\left\{l \left(\prod_{i=1}^{n-k-1} p_i\right) + x_0 : l<p_{n-k}\right\}$
  \item For all $x \in V_{\prod_{i=1}^{n-k}p_i}$, it exists a set of integers, $L$, such that for all $l \in L$ we have

$$\left(l \left(\prod_{i=1}^{n-k-1} p_i\right) + x_0\right)_{mod(\prod_{i=1}^{n-k}p_i)} \in U^{k+1}_{\prod_{i=1}^{n+1} p_i},$$

where $ \text{Card}(L) \cap [0:1:p_{n-k}-1]\geq p_{n-k}-k-1.$

\item The Pigeonhole Principle to conclude these facts if $k+1<\frac{p_{n-k}}{2}.$
\end{itemize} 

We would like to investigate if we can push the condition \ref{condition1} further, i.e., can we generalize our result for all $k<k_0$, for where $\prod_{i=1}^{n-k_0} p_i<p_{n+1}^2$ and $\prod_{i=1}^{n-k_0} p_i>p_{n+1}$ in order to ensure that set is not empty and contains only primes?\\

Our target is to push the condition \ref{condition1} more further, i.e., $k_0>\frac{p_{n-k}}{2}$, so we will prove that for all $k+m<k_0$ and for all
 $(x,y) \in \left(V_{\prod_{i=1}^{n-k} p_i}\right)^2$ then it exits $l \in \Z$ such that
\begin{equation}\label{target}
  \left\{x+l\prod_{i=1}^{n-k-m} p_i,\,y+l\prod_{i=1}^{n-k-m} p_i \right\}_{mod(\prod_{i=1}^{n-k} p_)} \subseteq \left(U^{k+1}_{\prod_{i=1}^{n+1} p_i}\right)^2
\end{equation}

In order to prove the previous equation some definitions and statement will be introduced to facilitate the explanation, then for all $(x,y) \in \left(V_{\prod_{i=1}^{n-k} p_i}\right)^2$ we define the following items

\begin{itemize}
  \item let $x_0=x_{mod(\prod_{i=1}^{n-k-m} p_i)}$ and $y_0=y_{mod(\prod_{i=1}^{n-k-m} p_i)}$
  \item let $V_{\prod_{i=1}^{n-k} p_i}=\left(V_{\prod_{i=1}^{n-k-m} p_i}+\prod_{i=1}^{n-k-m} p_i [0:1:\prod_{i=n-k-m+1}^{n-k} p_i]\right)/M_{\prod_{i=n-k-m+1}^{n+1}p_i},$ where $M_{\prod_{i=n-k-m+1}^{n+1}p_i}$ is a set of elements that are multiples of $p_{n-k-m+1}$ or $p_{n-k-m+2}$ or...or $p_{n+1}$
  \item let define the vector $H_x=(x_0, x_0+\prod_{i=1}^{n-k-m} p_i,...,x_o+l_x \prod_{i=1}^{n-k-m} p_i)$, where $l_x=\prod_{i=n-k-m+1}^{n-k} p_i-1$.
  \item let define the vector $H_y=(y_0, y_0+\prod_{i=1}^{n-k-m} p_i,...,y_o+l_y \prod_{i=1}^{n-k-m} p_i)$, where $l_y=\prod_{i=n-k-m+1}^{n-k} p_i-1$.
  \item For $j\in \{1,2,..,k+m \}$ and for $l \in \N$, if $\left(x_0+ l \prod_{i=1}^{n-k-m} p_i\right)_{mod(p_{n+1-j})}=0$ then $$\left(x_0+ (l\pm p_{n+1-j}) \prod_{i=1}^{n-k-m} p_i\right)_{mod(p_{n+1-j})}=0.$$
\end{itemize}

From the equation \ref{property2} and from the pigeonhole principle we can make ensure about the existence of $l \in \Z:\, x_0+l\prod_{i=1}^{n-k-m} p_i \in H_x,$ and  $y+l\prod_{i=1}^{n-k-m} p_i \in H_y$, e.i., the existence of $l' \in \Z$ such that $(x+l'\prod_{i=1}^{n-k-m} p_i,y+l'\prod_{i=1}^{n-k-m} p_i ) \in U^{k+1}_{\prod_{i=1}^{n+1} p_i} $ by choosing the right integers $m$ and $k$ such that

$$\text{Card} \left\{ l\in H_x \,:\, \forall i \in [n-k-m+1:1:n+1], gcd(l,p_i)=1 \right\}>\frac{\prod_{i=n-k-m+1}^{n-k} p_i}{2}.$$

By using Chinese remainder theorem, the density of co-primes, Lemma \ref{density_coprimes}, and $p_{n-k}\sim (n-k)\ln(n-k)$ to show that for $n-k$ big enough and for $m$ equals to two we have the following steps

\begin{description}
	\item[Step 1] if $\left(x_0+ l_j \prod_{i=1}^{n-k-2} p_i\right)_{mod(p_{n+1-j})}=0$ then for $i\in \{0,1,2,...,p_{n+1-j}\} $
	$$\left(x_0+ (l_j+i) \prod_{i=1}^{n-k-2} p_i\right)_{mod(p_{n+1-j})}\neq 0$$
	and
	$$\left(x_0+ (l_j+p_{n+1-j}) \prod_{i=1}^{n-k-2} p_i\right)_{mod(p_{n+1-j})}= 0.$$
	
	\item[Step 2] The existence of an integer number, $\Gamma$, such that
	\begin{eqnarray*}
	\left(\Gamma\right)_{mod(p_{n+1})} &=&p_{n+1}-\left(l_{n+1}\right)_{mod(p_{n+1})}\\
	\left(\Gamma\right)_{mod(p_{n})} &=& p_{n}-\left(l_{n}\right)_{mod(p_{n})}\\
	\left(\Gamma\right)_{mod(p_{n-1})} &=& p_{n-1}-\left(l_{n-1}\right)_{mod(p_{n-1})}\\
	&\vdots&\\
	\left(\Gamma\right)_{mod(p_{n-k-1})} &=& p_{n-k-1}-\left(l_{n-k-1}\right)_{mod(p_{n-k-1})}\\
	\end{eqnarray*}
is assured by the Chinese remainder Theorem.
\item[Step 3]  By consequence, for any positive integer $u$ we have the following
$$\text{Card}\left( l \in H'_{x} \right) \sim \text{Card}\left(\Gamma_u \leq j \leq\Gamma_u+\prod_{i=n-k-1}^{n-k} p_i\,|\, \forall i \in \{n-k-1,...,n+1\}\; gcd(j,p_i)=1 \right),$$
where $H'_{x}=\left\{l \in H_{x}\,|\,\forall i \in [n-k-1:1:n+1]\,:\, gcd(l,p_i)=1\right\}$ and $\Gamma_u=\Gamma +u \prod_{i=n-k-1}^{n+1} p_i.$
\item[Step 4] If $Q=\left\{p_i : n-k-1\leq i\leq n+1 \right\}$, then by using Lemma \ref{density_coprimes}  or the appendex we can write the following
\begin{eqnarray*}
  \text{Card} \left\{ l\in H'_x \right\} &=& \Psi(\Gamma_u+p_{n-k-1}p_{n-k},p_{n-k-2},Q) -\Psi(\Gamma_u+p_{n-k-2},p_{n-k-2},Q)\\
	&=& \Psi(p_{n-k-1}p_{n-k},p_{n-k-2},Q) -p_{n-k-2}\\
	&\sim&\left(p_{n-k-1}p_{n-k}-p_{n-k-2}\right)\left( \frac{\ln({n-k-2})}{\ln({n-k-1})(n-k)}\right)\\
   &\sim & \ln({n-k-2}) \frac{p_{n-k-1}p_{n-k}}{\ln(({n-k-1})({n-k}))} \\
   &\sim & \frac{p_{n-k-1}p_{n-k}}{2} \\
\end{eqnarray*}
\end{description}

Therefore,  it concludes for all $k\geq k_0$, where $k_0$ is smallest value that verify $\Pi_{i=1}^{n-k-2} p_i>p_{n+2}$, that

\begin{equation}\label{recursive_diff1}
  \left(\Delta V_{\prod_{i=1}^{n+1} p_i}\right)_{mod(\prod_{i=1}^{n+1} p_i)}=\left(\Delta V_{\prod_{i=1}^{n+1-k} p_i}\right)_{mod(\prod_{i=1}^{n} p_i)}+\sum_{j=0}^{k-1} \left(\prod_{i=1}^{n-j} p_i \right)[0:1:p_{n+1-j}-1]
\end{equation}

and
\begin{equation}\label{recursive_diff2}
\left(\Delta V_{\prod_{i=1}^{n+1} p_i}\right)_{mod(\prod_{i=1}^{n+1} p_i)}=\left(\Delta U^k_{\prod_{i=1}^{n+1} p_i}\right)_{mod(\prod_{i=1}^{n+1-k} p_i)}+\sum_{j=0}^{k-1} \left(\prod_{i=1}^{n-j} p_i \right)[0:1:p_{n+1-j}-1].
\end{equation}

From Lemma \ref{lemma 1} we can write
\begin{eqnarray*}
  \left(\Delta U^k_{\prod_{i=1}^{n+1} p_i}\right)_{mod(\prod_{i=1}^{n+1} p_i)}  &=& \left(\Delta U^{k,j}_{\prod_{i=1}^{n+1} p_i}\right)_{mod(\prod_{i=1}^{n+1} p_i)} \\
	&=&\left(\Delta V_{\prod_{i=1}^{n+1-k} p_i} \right)_{mod(\prod_{i=1}^{n+1-k}p_i)} \\
   &=& [0:2:\left(\prod_{i=1}^{n+1-k}p_i\right) -2].
\end{eqnarray*}
\newline

In order to omit the modulo from the previous equation, it will be enough to enlarge the set $U$ by a double. The calculation in modulo $\prod_{i=1}^{n} p_i$ will add $\prod_{i=1}^{n} p_i$ to the negative value, this idea will be used to prove the formula for the case $k=1$

$$ \left(\Delta U^1_{\prod_{i=1}^{n+1} p_i}\right) = \left(\Delta V_{\prod_{i=1}^{n+1} p_i} \right)_{mod(\prod_{i=1}^{n}p_i)},$$

We have the following facts:
\begin{itemize}
  \item $V_{\prod_{i=1}^{n+1} p_i} \bigcap [0:\prod_{i=1}^{n}p_i]=\bigcup_{j=1}^{p_{n+1}-1} \Psi_j$.
  \item $V_{\prod_{i=1}^{n+1} p_i} \bigcap [\prod_{i=1}^{n}p_i:2 \prod_{i=1}^{n}p_i]=(\Psi_0 +\prod_{i=1}^{n}p_i )\bigcup \left(\bigcup_{j=2}^{p_{n+1}-1} (\Psi_j+\prod_{i=1}^{n}p_i) \right)$.
  \item $\Psi_0$ is set of integers which are multiples of $p_{n+1}$ and belong to $V_{\prod_{i=1}^{n} p_i}$.
  \item In general $\Psi_j$ are sets of entegers which are multiples of $p_{n+1}$ and belong to $V_{\prod_{i=1}^{n} p_i}+j \prod_{i=1}^{n}p_i$, for $j=0,1,...,p_{n+1}-1$.
  \item If $x \in \Psi_0$ which means that $x_{mod(p_{n+1})}=0$ then $x+\prod_{i=1}^{n}p_i$ and $2\prod_{i=1}^{n}p_i -x$ both belong to $V_{\prod_{i=1}^{n+1} p_i} \bigcap [\prod_{i=1}^{n}p_i:2 \prod_{i=1}^{n}p_i].$
  \item If $x \in \Psi_1$ which means that $(x+\prod_{i=1}^{n}p_i)_{mod(p_{n+1})}=0$ then $x$ and $\prod_{i=1}^{n}p_i -x$ both belong to $V_{\prod_{i=1}^{n+1} p_i} \bigcap [0:\prod_{i=1}^{n}p_i].$
\end{itemize}

The question is how to make sure that for all $x$ and $y$ inside $\bigcup_{j=1}^{p_{n+1}-1} \Psi_j$ such that $x-y<0$ then $x+\prod_{i=1}^{n}p_i$ or $y-\prod_{i=1}^{n}p_i$ belong to $V_{\prod_{i=1}^{n+1} p_i} \bigcap [0:2 \prod_{i=1}^{n}p_i]$, i.e.,  $x-y+\prod_{i=1}^{n}p_i\in \Delta U^1_{\prod_{i=1}^{n+1}p_i}$.

To answer this question, we have to check the following cases

\begin{description}
  \item[Case 1] If $x$ and $y$ belong to $\bigcup_{j=2}^{p_{n+1}-1} \Psi_j$ then $x+\prod_{i=1}^{n}p_i$, $y+\prod_{i=1}^{n}p_i$ beleong to $V_{\prod_{i=1}^{n+1} p_i} \bigcap [0: \prod_{i=1}^{n}p_i]$. So $x-y+\prod_{i=1}^{n}p_i\in \Delta U^1_{\prod_{i=1}^{n+1}p_i}$.
  \item[Case 2] If  $x$ and $y$ belong to $\Psi_0$ then $\prod_{i=1}^{n}p_i-x$ and $\prod_{i=1}^{n}p_i-y$ belong to $\bigcup_{j=1}^{p_{n+1}-1} \Psi_j$ and also $x+\prod_{i=1}^{n}p_i$, $y+\prod_{i=1}^{n}p_i$, $2 \prod_{i=1}^{n}p_i-x$, and $2\prod_{i=1}^{n}p_i -y$ beleong to $V_{\prod_{i=1}^{n+1} p_i} \bigcap [\prod_{i=1}^{n}p_i:2 \prod_{i=1}^{n}p_i]$. So $x-y+\prod_{i=1}^{n}p_i\in \Delta U^1_{\prod_{i=1}^{n+1}p_i}$.
  \item[Case 3] If $x$ and $y$ belong to $\Psi_1$ then  $x$, $y$, $\prod_{i=1}^{n}p_i-x$, and $\prod_{i=1}^{n}p_i-y$ belong to $\bigcup_{j=1}^{p_{n+1}-1} \Psi_j$ and also  $2 \prod_{i=1}^{n}p_i-x$, and $2\prod_{i=1}^{n}p_i -y$ beleong to $V_{\prod_{i=1}^{n+1} p_i} \bigcap [\prod_{i=1}^{n}p_i:2 \prod_{i=1}^{n}p_i]$. So $x-y+\prod_{i=1}^{n}p_i\in \Delta U^1_{\prod_{i=1}^{n+1}p_i}$.
  \item[Case 4] If $x\in \Psi_1$ and $y \in \Psi_0 $ then  $x$, $x-\prod_{i=1}^{n}p_i$, and $\prod_{i=1}^{n}p_i-y$ belong to $\bigcup_{j=1}^{p_{n+1}-1} \Psi_j$ and also  $ 2\prod_{i=1}^{n}p_i-x$ and $y +\prod_{i=1}^{n}p_i$ beleong to $V_{\prod_{i=1}^{n+1} p_i} \bigcap [\prod_{i=1}^{n}p_i:2 \prod_{i=1}^{n}p_i]$. So $x-y+\prod_{i=1}^{n}p_i\in \Delta U^1_{\prod_{i=1}^{n+1}p_i}$.
  \item[Case 5] Other cases are trivial to be checked.
\end{description}\
\newline

For the other cases of $k$ can be derived by the same approach as we did to prove the equation \ref{target}, for $m=2$, to prove that

$ \forall x \in V_{\prod_{i=1}^{n-k} p_i}$ and $ \forall y \in \prod_{i=1}^{n-k}p_i+ V_{\prod_{i=1}^{n-k} p_i}$ then it exists $l \in \Z$ such that $(x+l\prod_{i=1}^{n-k-m} p_i )\in U^{k+1}_{\prod_{i=1}^{n+1} p_i}$ and $(y+l\prod_{i=1}^{n-k-m} p_i)\in U^{k+1,1}_{\prod_{i=1}^{n+1} p_i}.$\\

$\blacksquare$
\
\newline
\\

\section{Proof of Lemma \ref{density_coprimes}}
Let $\pi(x,n)$ be the co-prime-counting function that gives the number of co-primes with $p_i, i>n$, and are less or equal to $x$ and $N$ be the rank of biggest prime less than the value $x$, we can observe easily that
\begin{eqnarray*}
     \pi(x,1) &=& \pi(x) \\
     \pi(x,N) &=& x
  \end{eqnarray*}
	
Number of co-primes with $p_{N}$ and less than $x$ is calculated by subtracting all multiples of $p_N$ from $x$ integers, $$\pi(x,N-1)=x-\left\lfloor \frac{x}{p_N}\right\rfloor,$$ by the same way we can find number of co-primes with $p_N$ and $p_{N-1}$ can be calculated as
$$\pi(x,N-2)=x-\left\lfloor \frac{x}{p_N}\right\rfloor-\left\lfloor \frac{x}{p_{N-1}}\right\rfloor +\left\lfloor \frac{x}{p_N p_{N-1}}\right\rfloor,$$
where we simply this formula as
$$\pi(x,N-2)\sim x\left(1-\frac{1}{p_N}\right)\left(1-\frac{1}{p_{N-1}}\right). $$ 

By iteration, we can find $\pi(x,n)$ as follows
$$\pi(x,n)=x-\sum_{j=1}^N (-1)^j \left(\sum_{n+1\leq i_1 <i_2<...<i_j\leq N} \left\lfloor \frac{x}{\prod_{l=1}^{j}p_{i_l}}\right\rfloor\right)$$\
\newline
\\
We can generate it to
\begin{equation} \label{coprime1}
	\Psi(x,y,Q)=x-\sum_{j=1}^{N_Q} (-1)^j \left(\sum_{1\leq i_1 <i_2<...<i_j\leq N_Q} \left\lfloor \frac{x}{\prod_{l=1}^{j}q_{i_l}}\right\rfloor\right),
\end{equation}
where $q_.$ are primes inside the set $Q$ and $N_Q=Card(Q)$.\
\newline
\\
\subsection{Proof of the first part of Lemma \ref{density_coprimes}}
The first part of Lemma \ref{density_coprimes} can be conducted by the fact that for all positive integer $i$ and for all prime, $p \in Q$, we have
$$ \left(y+i\right)_{mod(p)}=\left(y+k P_Q+i\right)_{mod(p)},$$
then $$\Psi(x,y,Q)=Card \left\{i \in (y+k P_Q, x+k P_Q] \,|\,\, gcd(i,P_Q)=1 \right\}.$$\
\newline

\subsection{Proof of the second part of Lemma \ref{density_coprimes}}
The second part of Lemma \ref{density_coprimes} can be proved by using the the equation \ref{coprime1} and the fact that for all integer $x$
$$\left\lfloor \frac{x+k P_Q}{\prod_{l=1}^{j}q_{i_l}}\right\rfloor=\left\lfloor \frac{x}{\prod_{l=1}^{j}q_{i_l}}\right\rfloor + \frac{k P_Q}{\prod_{l=1}^{j}q_{i_l}}$$

as follows

\begin{eqnarray*}
\Psi(x+k P_Q,y,Q) &=& x+k P_Q-\sum_{j=1}^{N_Q} (-1)^j \left(\sum_{1\leq i_1 <i_2<...<i_j\leq N_Q} \left\lfloor \frac{x+k P_Q}{\prod_{l=1}^{j}q_{i_l}}\right\rfloor\right)\\
 &=& x+k P_Q-\sum_{j=1}^{N_Q} (-1)^j \left(\sum_{1\leq i_1 <i_2<...<i_j\leq N_Q} \left\lfloor \frac{x}{\prod_{l=1}^{j}q_{i_l}}\right\rfloor + \frac{kP_Q}{\prod_{l=1}^{j}q_{i_l}}\right)
\end{eqnarray*}

and
\begin{eqnarray*}
\Psi(y+k P_Q,y,Q) &=& y+k P_Q-\sum_{j=1}^{N_Q} (-1)^j \left(\sum_{1\leq i_1 <i_2<...<i_j\leq N_Q} \left\lfloor \frac{y+k P_Q}{\prod_{l=1}^{j}q_{i_l}}\right\rfloor\right)\\
&=& y+k P_Q-\sum_{j=1}^{N_Q} (-1)^j \left(\sum_{1\leq i_1 <i_2<...<i_j\leq N_Q} \left\lfloor \frac{y}{\prod_{l=1}^{j}q_{i_l}}\right\rfloor+ \frac{k P_Q}{\prod_{l=1}^{j}q_{i_l}}\right).
\end{eqnarray*}

By Subtracting the last two equations we can find
\begin{eqnarray*}
\Psi(x+k P_Q,y,Q)-\Psi(y+k P_Q,y,Q)&=&x-\sum_{j=1}^{N_Q} (-1)^j \left(\sum_{1\leq i_1 <i_2<...<i_j\leq N_Q} \left\lfloor \frac{x}{\prod_{l=1}^{j}q_{i_l}}\right\rfloor\right)-y\\
&=&\Psi(x,y,Q)-y,
\end{eqnarray*}\
\newline

where $\sum_{j=1}^{N_Q} (-1)^j \left(\sum_{1\leq i_1 <i_2<...<i_j\leq N_Q} \left\lfloor \frac{y}{\prod_{l=1}^{j}q_{i_l}}\right\rfloor\right)=0$ since $y<min(Q).$\
\newline
\\
We can find also that  if $x<P_Q$ and $gcd(x,P_Q)=1$ then
	$$ \Psi(x,y,Q)+\Psi(P_Q-x,y,Q)=\Psi(P_Q,y,Q)+1=P_Q \prod_{p_i \in Q} \left(1-\frac{1}{p_i}\right).$$
	
Therefore, by using same the equation \ref{coprime1} we can write the following equations 	
\begin{eqnarray*}
\Psi( P_Q-x,y,Q)+\Psi(y,y,Q)&=& P_Q-\sum_{j=1}^{N_Q} (-1)^j \left(\sum_{1\leq i_1 <i_2<...<i_j\leq N_Q} \left\lfloor \frac{P_Q-x}{\prod_{l=1}^{j}q_{i_l}}+ \frac{x}{\prod_{l=1}^{j}q_{i_l}}  \right\rfloor\right)\\
&=&-\sum_{j=1}^{N_Q} (-1)^j \left(\sum_{1\leq i_1 <i_2<...<i_j\leq N_Q}  -\left\lceil \frac{x}{\prod_{l=1}^{j}q_{i_l}}\right\rceil+ \left\lfloor \frac{x}{\prod_{l=1}^{j}q_{i_l}}\right\rfloor  \right)\\
& &\;+ P_Q- \sum_{j=1}^{N_Q} (-1)^j \left(\sum_{1\leq i_1 <i_2<...<i_j\leq N_Q}  \frac{P_Q}{\prod_{l=1}^{j}q_{i_l}} \right)\\
&=& P_Q- \sum_{j=1}^{N_Q} (-1)^j \left(\sum_{1\leq i_1 <i_2<...<i_j\leq N_Q}  \frac{P_Q}{\prod_{l=1}^{j}q_{i_l}} \right)+1\\
&=& P_Q \prod_{p_i \in Q} \left(1-\frac{1}{p_i}\right) +1
\end{eqnarray*}\
\newline	

\subsection{Proof of the third part of Lemma \ref{density_coprimes}}
The third part of Lemma \ref{density_coprimes} will be proved by improving Theorem 1 and the proposition 3 in \cite{coprimePSI} and \cite{coprimePSI2}, where it stated that: There is some $\delta$, $0\leq \delta <1$, and some $A>0$ such that

$$\sum_{\substack{p\leq x \; \\ \; p \in Q}} \ln(p)=\delta x + O(x \ln^{-A}(x)),$$
then
$$\Psi(x,y,Q)=\frac{e^{\gamma \delta}}{\Gamma (1-\delta)} x \prod_{\substack{y< p\leq x \; \\ \; p \in Q}} \left(1-\frac{1}{p}\right) \left(1+ O(D)\right),$$
where $D:=\frac{\ln(y)}{\ln(x)}+\ln^{-B}(y)$ and $B:=min(A,1)$.\
\newline

It is possible to choose the integer $k$ big enough in order to verify the following equality

$$\sum_{\substack{p\leq x \; \\ \; p \in Q}} \ln(p)= O((y+k P_Q) \ln^{-A}(y+k P_Q)),$$
where the value $\delta$ equals to zero and $A>1$.\
\newline
\\

Therefore, it implies that
$$\Psi(x+k P_Q,y,Q)= (x+k P_Q)\prod_{\substack{y< p\leq x \; \\ \; p \in Q}} \left(1-\frac{1}{p}\right) \left(1+ O(D_1)\right),$$
and 
$$\Psi(y+k P_Q,y,Q)= (y+k P_Q) \prod_{\substack{y< p\leq x \; \\ \; p \in Q}} \left(1-\frac{1}{p}\right) \left(1+ O(D_2)\right).$$

Then we can find for $k$ big enough that

\begin{eqnarray*}
\Psi(x+k P_Q,y,Q)-\Psi(y+k P_Q,y,Q)&=& (x-y)\prod_{\substack{y< p\leq x \; \\ \; p \in Q}} \left(1-\frac{1}{p}\right) \left(1+ O\left(\frac{1}{ln(y)}\right)\right)\\
&=&\Psi(x,y,Q)-y.
\end{eqnarray*}\
\newline
\\
Then it concludes that
$$\Psi(x,y,Q)=(x-y)\prod_{\substack{y< p\leq x \; \\ \; p \in Q}} \left(1-\frac{1}{p}\right) \left(1+ O\left(\frac{1}{ln(y)}\right)\right) + y.$$\
\newline

\subsection{Proof of the fourth part of Lemma \ref{density_coprimes}}
To Prove the last point of the Lemma \ref{density_coprimes}, it suffices to notice that $p_i \sim i\ln(i)$ in order to quantify the $\prod_{\substack{y< p\leq x \; \\ \; p \in Q}} \left(1-\frac{1}{p}\right) $ by simple integration as follows

\begin{eqnarray*}
\prod_{\substack{y< p\leq x \; \\ \; p \in Q}} \left(1-\frac{1}{p}\right)&=& \exp\left(\sum_{\substack{y< p\leq x \; \\ \; p \in Q}} \ln\left(1-\frac{1}{p}\right) \right)\\
&=&     \exp\left(\sum_{i=y/\ln(y)}^{x/\ln(x)} \ln\left(1-\frac{1}{i\ln(i)}\right) \right)\\
&\sim & \exp\left(\sum_{i=y/\ln(y)}^{x/\ln(x)} \left(-\frac{1}{i\ln(i)}\right) \right)\\ 
&\sim& \exp\left(\int_{y/\ln(y)}^{x/\ln(x)} \left(-\frac{1}{w\ln(w)}\right) dw \right)\\ 
&\sim& \frac{\ln(y)}{\ln(x)}.
\end{eqnarray*}\
\newline
\\
\begin{flushright}
$\blacksquare$
\end{flushright}

\
\newline
\\

\section{Proof of the conjecture of Goldbach} 
To prove the conjecture that said that all even number can be written as difference or sum of two prime numbers. Lemma \ref{LemmaP} will be the main idea to conclude the conjecture.
The question will be about the possibility to to find $k$ such that

\begin{equation}
p_{n+1}\leq \left(2 \prod_{i=1}^{n-k+1}p_i  \right) \leq p_{n+1}^{2}, \\
\end{equation}
in order to make sure that all integers inside the set $U_{\prod^{n-k+1}_{i=1} p_i}^{k,2}$ are primes.\\

For n big enough, we have that $p_n \sim n \ln(n)$ see \cite{paper2}, so
$$\prod_{i=1}^{n-k}p_i \sim \exp\left(\sum_{i\leq n-k} \ln(p_i)\right) \sim \exp\left(\sum_{i\leq n-k} \ln\left(i \ln(i)\right)\right),$$

then we can find the following approximation

\begin{eqnarray*}
  \prod_{i=1}^{n-k}p_i &\sim& \exp((n-k)\ln(n-k)) \\
 \end{eqnarray*}
This result can be calculated by using stirling formula.\
\newline 
In order to find $k$ as an increased and not bounded function of $n$, the following inequality needs to be solved

\begin{equation*}
(n+1) \ln(n+1)\leq  \exp((n-k)\ln(n-k))\leq (n+1)^{2} \ln^2(n+1).
\end{equation*}

In order to make sure that all values in $U_{\prod^{n-k+1}_{i=1} p_i}^k$ are prime numbers and their difference form all even numbers less than $\prod^{n-k+1}_{i=1} p_i$, the last inequality must be satisfied, then it is possible to choose
$$k\sim n-\frac{3\ln(n)}{2\ln(\ln(n))},$$

or in general as $n-f(n)$ where $f(n)^{f(n)}\sim n^{\alpha}$, for $\alpha \in (1,2)$,  in order to get

$$\exp((n-k)\ln(n-k)) \sim n^{\alpha}.$$\\

Till now we have proved that the difference of prime number forms all even number, therefore to prove that the sum of prime numbers form also all even number from the proof of the difference of prime number forms all even number it is enough to show that en mode of $\prod_{i=1}^{n}p_i$ we have

$$\Delta V_{\prod_{i=1}^{n}p_i} =V_{\prod_{i=1}^{n}p_i}+ V_{\prod_{i=1}^{n}p_i},$$

which it is trivial by using the fact that for all  $x \in V_{\prod_{i=1}^{n}p_i}$ implies that $\prod_{i=1}^{n} p_i-x \in V_{\prod_{i=1}^{n}p_i}$, i.e.,

\begin{eqnarray*}
  \left(\Delta V_{\prod_{i=1}^{n}p_i}\right)_{mod(\prod_{i=1}^{n}p_i)} &=& \left(\Delta V_{\prod_{i=1}^{n}p_i}+\prod_{i=1}^{n}p_i \right)_{mod(\prod_{i=1}^{n}p_i)}\\
   &=& \left( V_{\prod_{i=1}^{n}p_i}-V_{\prod_{i=1}^{n}p_i}+\prod_{i=1}^{n}p_i \right)_{mod(\prod_{i=1}^{n}p_i)} \\
   &=& \left( V_{\prod_{i=1}^{n}p_i}+V_{\prod_{i=1}^{n}p_i} \right)_{mod(\prod_{i=1}^{n}p_i)} . 
\end{eqnarray*}

\begin{flushright}
$\blacksquare$
\end{flushright}
\
\newline
\\
\section{Proof of the conjecture of twin prime}
In this section I will prove the conjecture by two different ways, the first prove it is just a sequences of  the lemma \ref{LemmaP} directly.The second one shows my first ideas when I started to find the proof of the twin number conjecture.

\subsection{First Approach}
The lemma \ref{LemmaP} stats that

$$\left(\Delta U^k_{\prod_{i=1}^{n+1} p_i}\right)\supseteq \left[0:2:\prod_{i=1}^{n-k+1} p_i-2\right],$$
where $ U^k_{\prod_{i=1}^{n+1} p_i} = V_{\prod_{i=1}^{n+1} p_i} \cap \left[ 0, 2 {\prod_{i=1}^{n-k+1} p_i}\right].$

Since the set $ U^k_{\prod_{i=1}^{n+1} p_i} $ contains only co-primes numbers with all the first $n+1$ primes numbers, the two smallest values inside the set  are $1$ and $p_{n+2}$. Therefore, to ensure that the set $ U^k_{\prod_{i=1}^{n+1}} $ contains only primes number, it is enough to choose $k$ around $n-\frac{\ln(n)}{\ln (\ln(n))}$ in order to ensure that the biggest co-prime inside $ U^k_{\prod_{i=1}^{n+1} }$ is less than $p_{n+1}^2$.

Since the boundaries of the set $ U^k_{\prod_{i=1}^{n+1}} $, without considering the value $1$, are going to infinity when $n$ goes to infinity, the existence of infinity twin number is verified\
\newline
\\
Since there is no prime between two primes that defer by 4, so the same proof can be conducted easily in order to prove the existence of infinitely couples of primes such that $p_{n+1}-p_{n}=4$, for every integer $k$. we can conclude further that $p_n=6k+1$ and $p_{n+1}=6k+5$ in modulo 6.

$$\liminf_{n\rightarrow \oo} \left( p_{n+1}-p_n\right)=2.$$\\
$\blacksquare$
\subsection{Second Approach}
The proof will use the both previous lemmas and the following concepts:\\

\textbf{Firstly}: $\lim_{n \rightarrow \oo} \pi(p_{n+1}^{2})-\pi(p_{n+1})=\oo$, where $\pi(x)$ is the prime counting function i.e., number of primes less than $x$. (see \cite{paper2} and \cite{paper3})\\

\textbf{Secondly}: The existence solution for the following equation
$$x \left(\prod_{i=1}^{n}p_i\right)+ y p_{n+1}=2,$$
where $x \in \Delta [1:p_{n+1}-1]$ and $y \in \Delta V_{\prod_{i=1}^{n}p_i}$.\\

For $n$ big enough, We need to prove the existence of twin primes $q$ and $q+2$ that belongs to
$$V_{\prod_{i=1}^{n}p_i}\bigcap \left[p_n, p_{n+1}^{2} \right[.$$
Referring to Theorem \ref{theorem 1}, $\forall  x \in V_{\prod_{i=1}^{n}p_i} $ and $x<p_{n+1}^{2}$ then $x$ is a prime number.

To prove the conjecture it is enough to show, for all $n$ big enough, that

$$ \exists (q,q+2) \in \left(V_{\prod_{i=1}^{n}p_i} \bigcap \left[p_n, p_{n+1}^{2}\right]\right)^2.$$

Suppose that $x_1$ and $x_2$ belong to $V_{\prod_{i=1}^{n+1}p_i}$ so it exists $(l_1,l_2) \in \left[1:p_{n+1}-1\right]^2$ and $(cp_1,cp_2) \in \left(V_{\prod_{i=1}^{n}p_i}\right)^2$ such that for all $(m_1, m_2) \in \Z^2$ we have

\begin{equation}
\begin{cases}
x_1=\left({\prod_{i=1}^{n}p_i}\right) l_1 + p_{n+1} cp_1 + m_1 \left({\prod_{i=1}^{n+1}p_i}\right)\\
x_2=\left({\prod_{i=1}^{n}p_i}\right) l_2 + p_{n+1} cp_2 + m_2 \left({\prod_{i=1}^{n+1}p_i}\right)
\end{cases}
\label{twinN}
\end{equation}

Abstracting both equations at the system \ref{twinN}, we get

$$x_2-x_1=\left({\prod_{i=1}^{n}p_i}\right) (l_2-l_1) + p_{n+1} (cp_2-cp_1) + (m_2-m_1) \left({\prod_{i=1}^{n+1}p_i}\right),$$
and if $x_1$ and $x_1$ are twin primes then
$$ 2+ (m_1-m_2) \left({\prod_{i=1}^{n+1}p_i}\right)=\left(V_{\prod_{i=1}^{n}p_i}\right) (l_2-l_1) + p_{n+1} (cp_2-cp_1)$$

By introducing the mod, we obtain
 \begin{equation}2=\left(  \left(V_{\prod_{i=1}^{n}p_i}\right) (l_2-l_1) + p_{n+1} (cp_2-cp_1)  \right)_{mod\left(V_{\prod_{i=1}^{n+1}p_i}\right)}.
\label{twinN2}
\end{equation}

All parameters in equation \ref{twinN2} need to be calculated and to show that both $x_1$ and $x_2$ are primes by showing their belong to the interval $[p_n; p_{n+1}^{2}[$.

Let denote $\Delta l= l_2-l_1$ and $\Delta cp=cp_2-cp_1$. From lemma \ref{lemma 1} we have
\begin{equation}
\begin{cases}
\Delta cp \in \left[0:2:{\prod_{i=1}^{n}p_i}-2\right]\\
\Delta l \in \left[2-p_{n+1}: p_{n+1}-2\right]
\end{cases}
\end{equation}

The equation \ref{twinN2} can be written as follows:

\begin{equation}2+m'{\prod_{i=1}^{n+1}p_i}=\left(  \left({\prod_{i=1}^{n}p_i}\right) \Delta l + p_{n+1} \Delta cp  \right).
\label{twinN3}
\end{equation}

By dividing the previous equation \ref{twinN3} by 2, we get the following

\begin{equation}1+m'{\prod_{i=3}^{n+1}p_i}=\left(  \left({\prod_{i=3}^{n}p_i}\right) \Delta l + p_{n+1} \frac{ \Delta cp}{2}  \right).
\label{twinN4}
\end{equation}

According to Diophantine theorem, the equation \ref{twinN4} has a particular solution, it means $\exists \Delta l^*$ and $\Delta cp^*$ such that

$$2+m'{\prod_{i=1}^{n+1}p_i}=\left(  \left({\prod_{i=1}^{n}p_i}\right) \Delta l^* + p_{n+1} \Delta cp^*  \right)$$

Then the general solution will be under the following form
\begin{equation}
\begin{cases}
\Delta cp =\Delta cp^*+k \prod_{i=1}^{n}p_i\\
\Delta l= \Delta l^* +k p_{n+1}
\end{cases}
\end{equation}
where $k\in \Z$.

We need to check if it exists $k$ such that
\begin{equation}
\begin{cases}
2-p_{n+1}\leq \Delta l \leq p_{n+1}-2\\
\text{and}\\
\Delta cp \in [0:2:\prod_{i=1}^{n}p_i-2].
\end{cases}
\end{equation}

It is obvious to see that for all $\Delta l^* \in Z$, it exists  $k$ such that $\Delta l^* +k p_{n+1} \in [2-p_{n+1} p_{n+1}-2],$ because if $(\Delta l^* )_{mod(p_{n+1})}=r$ then $k$ can be chosen accordingly
\begin{equation}
\begin{cases}
\text{if } r=p_{n+1}-1  \Rightarrow k= \left\lfloor -\frac{\Delta l^* }{p_{n+1}}\right\rfloor -1\\
\text{if } r \in [0,p_{n+1}-2]  \Rightarrow k= \left\lfloor -\frac{\Delta l^* }{p_{n+1}}\right\rfloor.
\end{cases}
\end{equation}

Still to check the possibility that 
$$\Delta cp \in [0:2:\prod_{i=1}^{n}p_i-2],$$ 
for that let start with the fist condition about $\Delta l$ that

$$2-p_{n+1} \leq\Delta l^* +k p_{n+1} \leq p_{n+1}-2,$$

By multiplication all parts by $\prod_{i=1}^{n}p_i$ we get
$$2\prod_{i=1}^{n}p_i- \prod_{i=1}^{n+1}p_i \leq \prod_{i=1}^{n}p_i\Delta l^* +k \prod_{i=1}^{n+1}p_i \leq \prod_{i=1}^{n+1}p_i-2\prod_{i=1}^{n}p_i.$$

Indeed, it follows from equation \ref{twinN3} that $2=\left(  \left({\prod_{i=1}^{n}p_i}\right) \Delta l^* + p_{n+1} \Delta cp^*  \right)- k_0 (\prod_{i=1}^{n+1}p_i)$
then
$$2\prod_{i=1}^{n}p_i- \prod_{i=1}^{n+1}p_i \leq  \left(2-p_{n+1} \Delta cp^*\right) +(k+k_0) \prod_{i=1}^{n+1}p_i \leq \prod_{i=1}^{n+1}p_i-2\prod_{i=1}^{n}p_i,$$
therefore
$$\frac{2\prod_{i=1}^{n}p_i -2}{p_{n+1}}- \prod_{i=1}^{n}p_i \leq -\Delta cp^* +(k+k_0) \prod_{i=1}^{n}p_i \leq \frac{-2\prod_{i=1}^{n}p_i +2}{p_{n+1}}+ \prod_{i=1}^{n}p_i.$$

For all even possible value of $\Delta cp^*$, we can find $k_0$ in order to make sure that $\Delta cp \in [0:2:\prod_{i=1}^n p_i-2]$.

Finally $x_1$ and $x_2$ can be written as

\begin{equation}
\begin{cases}
x_1=\left({\prod_{i=1}^{n}p_i}\right) l_1 + p_{n+1} cp_1 + m_1 \left({\prod_{i=1}^{n+1}p_i}\right)\\
x_2=\left({\prod_{i=1}^{n}p_i}\right) \left(\Delta l^* +l_1\right) + p_{n+1} \left(\Delta cp^* + cp_1\right) + m_2 \left({\prod_{i=1}^{n+1}p_i}\right)
\end{cases}
\label{twinNs}
\end{equation}

Another question is raised to check and verify again the following
\begin{equation}
\begin{cases}
\left(\Delta l^* +l_1\right)_{mod(p_{n+1})} \in [1:p_{n+1}-1]\\
\text{and }\\
\left(\Delta cp^* + cp_1\right)_{mod(\prod_{i=1}^{n}p_{i})} \in V_{mod(\prod_{i=1}^{n}p_{i})}\\
\end{cases}
\label{twinNs2}
\end{equation}

From Theorem \ref{theorem 1} , we need to be sure that $\Delta l^* +l_1$ and $p_{n+1}$ are co-primes and also $\Delta cp^* + cp_1$ is co-prime with $\prod_{i=1}^{n}p_{i}$. To do that, we can see that $\Delta l^* $ is co-prime with $p_{n+1}$ and  $\Delta cp^*/2 $ is co-prime with $\prod_{i=3}^{n}p_{i}$, then $l_1$ can be chosen easily in order to avoid that $\left(\Delta l^* +l_1\right)_{mod(p_{n+1})}=0$.\\
\
\newline

The only question of the existence of $cp_1$ such that   $cp_1 \in V_{\prod_{i=3}^{n}p_{i}}$, $\Delta cp^* + cp_1 \in V_{\prod_{i=1}^{n}p_{i}}$, and $x_1 \in (p_{n+1}, p_{n+1}^2)$ needs to be answer in order to finish the prove of the conjecture.

It was proved that gap between primes is finite,\cite{paper4}, by consequence  for $n$ big enough, $n \sim p_{\left[\Pi(n)\right]}$, then by substitution  the variable $n$ by  $n \ln(n)$, we get very interesting formula $$p_n\sim n \ln(n).$$

Number of prime numbers less than $n$ can be calculated differently by

$$\pi(n)=n-card\left(\bigcup_{i=2}^n A_i\right)=n-\left(\sum_{k=2}^{n} (-1)^{k}  card\left( \bigcap_{i=2}^{k} A_i\right) \right),$$

where $A_i=\left\{{m\,\, :\, m\in [i+1,n] \text{ and } m_{mod(i)}=0}\right\}$, in our next paper, a new method in how to estimate $\pi(n)$ will be our focus, here are some headline

\begin{eqnarray*}
\pi(n)&=& n-\left(\sum_{k=2}^{n} (-1)^{k}  card\left( \bigcap_{i=2}^{k} A_i\right) \right)\\
&=& n-\left(\sum_{k=2}^{n} (-1)^{k} \sum_{2\leq i_1<i_2<...<i_k\leq n} \left(card\left( \left\lfloor \frac{n}{\prod_{i=2}^{k} n_i} \right\rfloor -1\right) \right)\right)\\
&=&  n-\left(\sum_{k=2}^{n} (-1)^{k} \sum_{2\leq i_1<i_2<...<i_k\leq n} \left(card\left( \left\lfloor \frac{n}{\prod_{i=2}^{k} n_i} \right\rfloor \right) \right)\right) -1\\
&\approx & n \prod_{i=2}^{n} \left(1-\frac{1}{i}\right)\\
&\approx & \frac{n}{\ln(n)}.
\end{eqnarray*}\
\newline

Furthermore, the existence of several choices of $cp_1$ such that $\Delta cp^* + cp_1 \in V_{\prod_{i=1}^{n}p_{i}}$ is assured by Lemma \ref{LemmaP}. The value of $\Delta cp^*$ will be successively divided by $\prod_{i=1}^{n}p_{i}$ as follows:

\begin{eqnarray*}
  \Delta cp^* &=& l_1 \prod_{i=1}^{n-1}p_{i} +r_1 \\
  r_1 &=& l_2 \prod_{i=1}^{n-2}p_{i} +r_2 \\
   & \vdots &  \\
  r_{k-1} &=& l_k \prod_{i=1}^{n-k}p_{i} +r_{k}
\end{eqnarray*}

where $r_.$ and $l_.$ are the quotients and reminders of each division.

Then $\Delta cp^*$ can be written as follows:

$$\Delta cp^* =r_k+ \sum_{j=1}^{k} l_j \left(\prod_{i=1}^{n-j}p_{i}\right) $$

The value of $k$ should be chosen in order to be sure that $x_1$ belongs to the interval of primes number only, $[p_{n+1}, p^2_{n+1}) \cap V_{\prod_{i=1}^{n+1}p_{i}}$. To do so $k$ needs to be the first value that verified the following inequality:

$$\left(\prod_{i=1}^{n-k}p_{i}\right)\leq \frac{p_{n+1}}{2}.$$

From the part b of Lemma \ref{lemma 1}, it is clear that it exists an integer $j$, strictly positive, for which
$$cp_1 \in \left[\frac{{j \prod_{i=1}^{n}p_{i}}}{p_{n+1}} , \frac{{j \prod_{i=1}^{n}p_{i}}}{p_{n+1}}+ p_{n+1}\right) \cap V_{\prod_{i=1}^{n+1}p_{i}},$$
and another co-prime number, $cp_2 \in V_{\prod_{i=1}^{n+1}p_{i}}$ such that $\Delta cp^*=cp_2-cp_1$

We are sure now about the possibility of selecting the right $l_1$ in order to make sure that $x_1$ is a prime number by assuring that $x_1$ is inside the interval $[p_{n+1}, p^2_{n+2}) \cap V_{\prod_{i=1}^{n+1}p_{i}}$.\\

Since there is no prime between two primes that defer by 4, so the same proof can be conducted easily in order to prove the existence of infinitely couples of primes such that $p_{n+1}-p_{n}=4$, for every integer $k$. we can conclude also that $p_n=1$ in modulo 6.

$\blacksquare$

\
\newline
\\
\section{Appendex}
In this section, new approach in how to calculate the density of co-primes will be described. Let's suppose that for all integer $x>p_n$ we have for $Q_m=\left\{p_i\,|\, n\leq i< n+m \right\}$ and $p_N$ represent the biggest prime less than $x$

$$\Psi(x,p_{n-1},Q_m)=x \prod_{i=n}^{n+m-1} \left(1-\frac{1}{p_i}\right)+ O(f(x,m)), \text{ for all } m \leq N-n. $$ 

Using the proof by induction, the function $f(x,m)$ can be estimated as a linear function with $\ln(m+n)$.\\

For $m=1$, we have

\begin{eqnarray*}
     \Psi(x,p_{n-1},Q_1) &=& x \prod_{i=n}^{n} \left(1-\frac{1}{p_i}\right)+ O(f(x,1))\\
     &=& x-\frac{x}{p_n}+O(f(x,1)),
  \end{eqnarray*}

The error function, $O(f(x,1))$, can be calculated easily as follows

\begin{eqnarray*}
     O(f(x,1)) &=& x-\left\lfloor \frac{x}{p_n}\right\rfloor-\left( x-\frac{x}{p_n}\right) \\
     &=& \frac{x}{p_n}- \left\lfloor \frac{x}{p_n}\right\rfloor\\
		&\leq& 1.
  \end{eqnarray*}\
	\newline
	
Let's assume that $$\Psi(x,p_{n-1},Q_m)=x \prod_{i=n}^{n+m-1} \left(1-\frac{1}{p_i}\right)+ O(f(x,m))$$
and  we need to show

$$\Psi(x,p_{n-1},Q_{m+1})=x \prod_{i=n}^{n+m} \left(1-\frac{1}{p_i}\right)+ O(f(x,m+1)).$$

Let $A_i$ be set of multiple of $p_{n+i}$ and less than $x$, $A_i=\left\{i\,|\, gcd(i,p_{n+i-1})=p_{n+i-1} \, \text{ and } \, i\leq x \right\}$, then

\begin{eqnarray*}
     \Psi(x,p_{n-1},Q_{m+1}) &=&  x- card\left(\bigcup_{i=1}^{m+1} A_i\right)\\
     &=& x- \text{card}\left(\left(\bigcup_{i=1}^{m} A_i \right) \bigcup  A_{m+1}\right)\\
		&=& x- \text{card}\left(\left(\bigcup_{i=1}^{m} A_i \right)\right) - \text{card}\left( A_{m+1} \right) + \text{card}\left(\left(\bigcup_{i=1}^{m} A_i \right) \bigcap  A_{m+1}\right)\\
		&=& x- \text{card}\left(\left(\bigcup_{i=1}^{m} A_i \right)\right) - \left\lfloor \frac{x}{p_{n+m}}\right\rfloor + \text{card}\left(\left(\bigcup_{i=1}^{m} A_i \cap A_{m+1}\right)   \right).
  \end{eqnarray*}\
	\newline
Note that $A_i \cap A_{m+1}=\left\{i\,|\, gcd(i,p_{n+m} p_{n+i-1})=p_{n+m }p_{n+i-1} \, \text{ and } \, i\leq x \right\}$ and 
$\Psi(x,p_{n-1},Q_{m})=x- \text{card}\left(\left(\bigcup_{i=1}^{m} A_i \right)\right)$ then

\begin{eqnarray*}
     \Psi(x,p_{n-1},Q_{m+1}) &=& \Psi(x,p_{n-1},Q_{m}) - \left\lfloor \frac{x}{p_{n+m}}\right\rfloor +\frac{x}{p_{n+m}}- \Psi\left(\frac{x}{p_{n+m}},p_{n-1},Q_{m}\right) \\
		&=& x \prod_{i=n}^{n+m-1} \left(1-\frac{1}{p_i}\right) - \left\lfloor \frac{x}{p_{n+m}}\right\rfloor +\frac{x}{p_{n+m}}- \frac{x}{p_{n+m}} \prod_{i=n}^{n+m-1} \left(1-\frac{1}{p_i}\right)\\
		& &\; + O \left(f(x,m)-f(x/p_{n+m},m)\right)\\
		&=& x \prod_{i=n}^{n+m} \left(1-\frac{1}{p_i}\right) + O \left(f(x,m)-f(x/p_{n+m},m)+1\right)\\
  \end{eqnarray*}\
	\newline
In order to force the last formula to correct we have to find the function, $f$, that verify the following

$$ f(x,m)+f(x/p_{n+m},m)=f(x,m+1) +O(1).$$

For $n$ is big enough integer, then
\begin{eqnarray*}
    -f(x/p_{n+m},m)  &=& f(x,m+1)-f(x,m) \\
		& \sim & \frac{\partial f(x,m)}{\partial m}
  \end{eqnarray*}\
	\newline
	
Since $p_{n+m}\sim (n+m)\ln(n+m)$ then differential equation needs to be solved

$$\frac{\partial f(x,m)}{\partial m}= -f\left(\frac{x}{(n+m)\ln(n+m)},m\right).$$

The function $$f(x,m) \sim \frac{x}{\ln(x)\, \ln(n+m)}$$ is a possible because

\begin{eqnarray*}
\frac{x}{\ln(x)\, \ln(n+m)} -\frac{x/p_{n+m}}{\ln(x/p_{n+m})\, \ln(n+m)}  &=& \frac{x}{\ln(x)} \left( \frac{1}{\ln(n+m)}\right)\\
& & -\left( \frac{x}{(n+m)\ln(x) \ln^2(n+m) \left(1+\ln(n+m)/\ln(x)\right)}\right) \\
&\leq& \frac{x}{\ln(x)} \left( \frac{1}{\ln(n+m)}-\frac{1}{2(n+m)\ln^2(n+m) }\right)\\
		&\leq& \frac{x}{c \ln(x)\, \ln(n+m+1}\\
		&\leq& f(x,m+1)
  \end{eqnarray*}
	
	where $c$ is a fixed constant bigger than one.
%
	%
	%
then
$$\Psi(x,p_{n-1},Q_m)=x \prod_{i=n}^{n+m-1} \left(1-\frac{1}{p_i}\right)+ O\left(\frac{x}{\ln(x)\, \ln(n+m)}\right), \text{ for all } m \leq N-n. $$

For $n$ big enough we can write
$$ \Psi(x,p_{n-1},Q_{N-n})=x \frac{\ln(n)}{\ln(x)}+ O\left(\frac{x}{\ln^2(x)}\right).$$
\begin{flushright}
$\blacksquare$
\end{flushright}\
\newline
\\

\end{document}